\documentclass[reqno,a4paper]{amsart}
\usepackage{graphicx}
\usepackage{amsmath}
\usepackage{mathrsfs}
\usepackage{amsfonts}
\usepackage{amssymb}
\usepackage{enumerate}
\usepackage{latexsym}
\usepackage{graphpap}
\usepackage{cite}
\usepackage{tikz}

\allowdisplaybreaks
\newtheorem{theorem}{Theorem}[section]

\newtheorem{lemma}[theorem]{Lemma}
\newtheorem{proposition}[theorem]{Proposition}

\theoremstyle{definition}

\newtheorem{example}[theorem]{Example}

\newtheorem{remark}[theorem]{Remark}

\newtheorem{algo}{Algorithm}

\theoremstyle{remark}
\numberwithin{equation}{section}

\newcommand{\bfa}{\mathbf{a}}
\newcommand{\bfb}{\mathbf{b}}

\newcommand{\bfp}{\mathbf{p}}
\newcommand{\bfq}{\mathbf{q}}

\newcommand{\bfs}{\mathbf{s}}
\newcommand{\bft}{\mathbf{t}}
\newcommand{\bfu}{\mathbf{u}}
\newcommand{\bfv}{\mathbf{v}}
\newcommand{\bfw}{\mathbf{w}}

\newcommand{\bfx}{\mathbf{x}}
\newcommand{\bfy}{\mathbf{y}}

\newcommand*{\hypo}{{\mathsf{hypo}}}

\newcommand{\ev}{\operatorname{\mathsf{ev}}}
\newcommand{\supp}{\operatorname{\mathsf{sup}}}

\newcommand{\rpi}{\operatorname{\mathsf{rpi}}}
\newcommand{\lpi}{\operatorname{\mathsf{lpi}}}
\newcommand{\fp}{\operatorname{\mathsf{fp}}}
\newcommand{\ip}{\operatorname{\mathsf{ip}}}

\newcommand{\op}{\,^*\,}
\newcommand{\var}{\mathsf{Var}}
\newcommand{\fb}{finitely based}
\newcommand{\nfb}{non-finitely based}
\newcommand{\inmon}{involution monoid}
\newcommand{\insem}{involution semigroup}

\newcommand{\con}{\operatorname{\mathsf{con}}}

\newcommand{\occ}{\operatorname{\mathsf{occ}}}

\newcommand{\sylv}{\operatorname{\mathsf{sylv}}}
\newcommand{\baxt}{\operatorname{\mathsf{baxt}}}

\begin{document}
\begin{sloppypar}
\title[Representations and identities of Baxter monoids with involution]{Representations and identities of Baxter monoids with involution}%
\thanks{The authors are partially supported by the National Natural Science Foundation of China (Nos. 12271224, 12171213, 12161062).}

\author[B. B. Han]{Bin Bin Han$^{1}$}
\author[W. T. Zhang]{Wen Ting Zhang$^{1,\star}$}\thanks{$^\star$ Corresponding author} %
\author[Y. F. Luo]{Yan Feng Luo$^{1}$}
\author[J. X. Zhao]{Jin Xing Zhao$^{2}$}

\address{$^{1}$ School of Mathematics and Statistics, Lanzhou University, Lanzhou, Gansu 730000, PR China} %
\email{hanbb19@lzu.edu.cn, zhangwt@lzu.edu.cn, luoyf@lzu.edu.cn}

\address{$^{2}$ School of Mathematical Sciences, Inner Mongolia University, Hohhot, Inner Mongolia 010021, PR China}.
\email{zhjxing@imu.edu.cn}



\begin{abstract}
Let $(\mathsf{baxt}_n,~^\sharp)$ be the Baxter monoid of finite rank $n$ with Sch\"{u}tzenberger's involution $^{\sharp}$.
In this paper, it is shown that $(\mathsf{baxt}_n,~^\sharp)$ admits a faithful representation by an involution monoid of upper triangular matrices over any semiring from a large class including the tropical semiring under the skew transposition. Then a transparent combinatorial characterization of the word identities satisfied by $(\mathsf{baxt}_n,~^\sharp)$ is given. Further, it is proved that $(\mathsf{baxt}_n,~^\sharp)$ is finitely based if and only if $n\neq 3$, and shown that the identity checking problem for $(\mathsf{baxt}_n,~^\sharp)$ can be done in polynomial time.
\end{abstract}
%\date{\today}
\maketitle

\section{Introduction}
Identities and varieties of semigroups have long been studied, and several important questions arise in this area, such as the \textit{finite basis problem}, that is the problem of classifying semigroups according to the finite basis property. The first example of {\nfb} finite semigroup was discovered by Perkins~\cite{Perkins69} in the 1960s; since then, the finite basis problem for semigroups has attracted much attention. Refer to Volkov~\cite{Vol01} for a survey of work performed in this direction and for more information on the finite basis problem for semigroups in general. Other questions regarding the variety generated by a semigroup are those of whether it contains only finitely generated subvarieties (see, for example, \cite{Vol01}), or countably infinite subvarieties \cite{Jack00}. Another widely studied question is the \textit{identity checking problem} \textsc{Check-Id}($S$) for a semigroup $S$, that is the decision problem whose
instance is an arbitrary identity $\bfu \approx \bfv$, and the answer to such an instance is `YES'
if $S$ satisfies $\bfu \approx \bfv$, and `NO' if it does not.
For a finite semigroup $S$, the identity checking problem \textsc{Check-Id}($S$) is always decidable  but this is not necessarily true for an infinite semigroup \cite{Mu68}.
Studying the computational complexity of identity checking in semigroups was proposed
by Sapir \cite[Problem~2.4]{KS95}, and up to now, there are many results in this study \cite{AVG09,JM06,DJK18,Chen2020, KV20}.

Recall that a unary operation~$^*$ on a semigroup~$S$ is an \textit{involution} if~$S$ satisfies the identities
\begin{equation}
(x^*)^* \approx  x \quad \text{and} \quad (xy)^* \approx y^*x^*. \label{id: inv}
\end{equation}
An \textit{involution semigroup} is a pair $(S, \op)$ where~$S$ is a semigroup with {involution}~$^*$, and
$S$ is called the \textit{semigroup reduct} of $(S,\op)$.  Common examples of involution semigroups include groups $(G,\,^{-1}\,)$ with inversion~$^{-1}$, multiplicative matrix semigroups $(M_n,\,^T\,)$ and $(M_n,\,^D\,)$ over any field with transposition~$^T$ and skew transposition~$^D$ respectively.
Over the years, the identities and varieties of {\insem}s have received less attention than that for semigroups.
Since the turn of the millennium, interest in {\insem}s has significantly increased.
For example, many counterintuitive results were established and examples have been found to demonstrate that an {\insem} $(S,\op)$ and its semigroup reduct~$S$ need not be simultaneously {\fb} \cite{GZL-TB2,GZL-A01,JV10,Lee16a,Lee16b,Lee17a,Lee18,Lee19}. Refer to \cite{Lee20,Lee19+,GZL-variety,ADPV14,ADV12a,ADV12b} for more information on the identities and varieties of involution semigroups.

The plactic monoid, whose elements can be identified with Young
tableaux, is famous for its interesting connection with combinatorics \cite{Lot02} and applications in symmetric functions \cite{Mac08}, representation theory \cite{Ful97}, Kostka-Foulkes polynomials \cite{LS78,LS81} and Schubert polynomials \cite{LS85, LS90}. Its finite-rank versions were shown to have faithful tropical representations by Johnson and Kambites \cite[Theorem~2.8]{JK19}. An important consequence of these representations, which are specifically representations using upper triangular tropical matrices, is that each finite-rank plactic monoid satisfies a non-trivial semigroup identity \cite[Theorem~3.1]{JK19}. In particular, the plactic monoids of rank equal to $2$ and $3$ are non-finitely based by \cite{CHLS16,DJK18} and \cite{HZL,JK19} respectively, and their involution cases are also non-finitely based by using the result of \cite[Theorem~4]{Lee17a}. The finite basis problem for the plactic monoids of rank greater than $3$ and their involution cases are still open. The plactic monoid of infinite rank is finitely based since it only satisfies trivial identity \cite[Proposition~3.1]{CKK17}.

Plactic monoids belong to a family of `plactic-like' monoids which are connected with
combinatorics and whose elements can be identified with combinatorial objects. Others in the family include the hypoplactic monoids, whose elements are quasi-ribbon tableaux and which play a role in the theory of quasi-symmetric functions analogous to that of the plactic monoid for symmetric functions \cite{KT97,Nov00}; the sylvester and ${\sharp}$-sylvester monoids, whose elements are respectively right strict and left strict binary search trees \cite{HNT05}; the Baxter monoids, whose elements are pairs of twin binary search trees \cite{Gir12} and which are linked to the theory of Baxter permutations; the taiga monoids, whose elements are binary search trees with multiplicities \cite{Pri13}; the stalactic monoids, whose elements are stalactic tableaux \cite{HNT07,Pri13}; and the stylic monoids, whose elements are $N$-tableaux \cite{AR22}.

The hypoplactic, stalactic, taiga, sylvester and Baxter monoids of each finite rank were shown to have faithful tropical representations \cite{CJKM21}, and each of these monoids satisfy non-trivial identities \cite{CM18}.
It is shown that all Baxter [resp. hypoplactic, sylvester, stalactic and taiga] monoids of rank greater than or equal to $2$ generate the same variety and are finitely based \cite{CMR21b, CJKM21,HZ,CMR21a}.
Aird and Ribeiro have given a faithful representation of the stylic monoid and its involution case of each finite rank, and then solved the finite basis problems for them \cite{TD22}.
Also, Volkov has solved the finite basis problem for the stylic monoid by different mean \cite{volkov2022}.
The identity checking problems for the Baxter, sylvester and stylic monoids have been considered \cite{CMR21b,TD22}, which can be done in polynomial time.
By the characterization of the identities satisfied by the hypoplactic, stalactic and taiga monoids \cite{CMR21b, CJKM21,HZ}, it is easy to see that the identity checking problem for each of them is in the complexity class $\mathsf{P}$.

The Baxter monoid $\mathsf{baxt}_n$ of finite rank $n$, first studied  by  Giraudo \cite{Gir12}, can be obtained by factoring the free monoid $\mathcal{A}_n^{\star}=\mathcal{A}_n^{+}\cup \{\varepsilon\}$ over the finite ordered alphabet $\mathcal{A}_n = \{1 < 2 < \cdots < n\}$ by a congruence that can be defined by the insertion algorithm that computes a combinatorial object from a word. Each element of the Baxter monoid can be uniquely identified with a pair of twin binary trees, and it is presented by the relations defined for $\bfu, \bfv \in \mathcal{A}_n^{\star}$ and $a, b, c, d\in \mathcal{A}_n$ by
$c\bfu ad\bfv b \equiv  c\bfu da\bfv b$ with $a \leq b < c \leq d$ and
$b\bfu da\bfv c \equiv b\bfu ad\bfv c$ with $a < b \leq c < d$. The Baxter monoid $\mathsf{baxt}_n$ forms the involution monoid $(\mathsf{baxt}_n, ~^{\sharp})$ under the order-reversing permutation on $\mathcal{A}_n$.

In this paper, we investigate the representations and identities of finite-rank Baxter monoid with involution $(\mathsf{baxt}_n,~^\sharp)$.
This paper is organized as follows. Notation and background information of the paper
are given in Section \ref{sec: prelim}. In Section \ref{sec:repre}, we exhibit a faithful representation of $(\mathsf{baxt}_n,~^\sharp)$ as an involution monoid of upper triangular matrices over any semiring
from a large class including the tropical semiring under the skew transposition. In Section \ref{sec:characterization},  a characterization of  the word identities satisfied by $(\mathsf{baxt}_n,~^\sharp)$ is given. In Section \ref{sec:FBP4},  we prove that $(\mathsf{baxt}_n,~^\sharp)$ is finitely based if and only if $n\neq3$  and that each variety generated by $(\mathsf{baxt}_n,~^\sharp)$ with $n\geq 2$ contains continuum many subvarieties. And it is shown that the identity checking problem for $(\mathsf{baxt}_n,~^\sharp)$ can be done in polynomial time in Section \ref{sec:CHECK-ID}.

\section{Preliminaries} \label{sec: prelim}
Most of the notation and definition of this article are given in this section.
Refer to the monograph of Burris and Sankappanavar~\cite{BS81} for any undefined notation and terminology of universal algebra in general.

\subsection{Words and content}

Let~$\mathcal{X}$ be a countably infinite alphabet and $\mathcal{X}^* =\{x^* \,|\, x \in \mathcal{X}\}$ be a disjoint copy of~$\mathcal{X}$.
Elements of $\mathcal{X} \cup \mathcal{X}^*$ are called \textit{variables}, elements of the free involution monoid $(\mathcal{X} \cup \mathcal{X}^*)^{\times}=(\mathcal{X} \cup \mathcal{X}^*)^+ \cup \{\emptyset\}$ are called \textit{words}, and elements of $\mathcal{X}^+ \cup \{\emptyset\}$ are called \textit{plain words}. A word $\bfu$ is a \textit{factor} of a word $\bfw$ if $\bfw = \bfa\bfu\bfb$ for some $\bfa, \bfb \in (\mathcal{X} \cup \mathcal{X}^*)^{\times}$.

Let $\bfu\in (\mathcal{X} \cup \mathcal{X}^*)^+$ be a word and $x \in \mathcal{X} \cup \mathcal{X}^*$ be a variable. The \textit{content} $\con (\bfu)$ of a word $\bfu$ is the set of variables that occur in~$\bfu$, and $\occ(x, \bfu)$ is the number of occurrences of the letter $x$ in $\bfu$. Let $\overline{\bfu}$ be the plain word obtained from $\bfu$ by removing all occurrences of the symbol $^*$.
For any $x_1, x_2, \ldots, x_n \in \mathcal{X}\cup\mathcal{X}^*$ such that $\overline{x_1}, \overline{x_2}, \ldots, \overline{x_n}\in \mathcal{X}$ are distinct variables, let $\bfu[x_1, x_2, \ldots, x_n]$ denote the word obtained from~$\bfu$ by retaining only the variables $x_1, x_1^*, x_2, x_2^*, \ldots, x_n, x_n^*$. Denote by $\overleftarrow{\occ}_y(x, \bfu)$ [resp. $\overrightarrow{\occ}_y (x, \bfu)$] the number of occurrences of $x$ appearing before [resp. after] the first [resp. last] occurrence of $y$ in $\bfu$.  The \textit{initial part} of $\bfu$, denoted by $\ip(\bfu)$, is the word obtained from $\bfu$ by retaining the occurrence of each variable $x\in\con(\bfu)$ satisfying $x, x^{* }\not\in\con(\bfu_1)$ where $\bfu=\bfu_1x\bfu_2$; the \textit{final part} of $\bfu$, denoted by $\fp(\bfu)$, is the word obtained from $\bfu$ by retaining the occurrence of each variable $x\in\con(\bfu)$ satisfying $x, x^{* }\not\in\con(\bfu_2)$ where $\bfu=\bfu_1x\bfu_2$.

\begin{example}
Let $\bfu = x^* zxy^*xyz^2x$ and $x,y,z\in \mathcal{X}$. Then
\begin{itemize}
  \item  $\con (\bfu)=\{x, y, z, x^*, y^*\}$, $\overline{\bfu}=xzxyxyz^2x$,
  \item $\occ (x, \bfu) = 3$, $\occ (x^*,\bfu) = \occ (y,\bfu)=\occ (y^*,\bfu)=1$,
  \item $\bfu[x]=x^*x^3$, $\bfu[x, y]=x^*xy^*xyx$,
  \item $\overleftarrow{\occ}_{y^*}(x, \bfu)=1$, $\overrightarrow{\occ}_{y^*}(x, \bfu)=2$,
  \item $\ip(\bfu)=x^*zy^*$, $\fp(\bfu)=yzx$.
\end{itemize}
\end{example}

\subsection{Terms and identities}

The set $\mathsf{T}(\mathcal{X})$ of \textit{terms} over~$\mathcal{X}$ is the smallest set containing $\mathcal{X}$ that is closed under concatenation and $^*$.
The proper inclusion $(\mathcal{X} \cup \mathcal{X}^*)^{\times} \subset \mathsf{T}(\mathcal{X})$ holds and the identities \eqref{id: inv} can be used to convert any nonempty term $\bft \in \mathsf{T}(\mathcal{X})$ into some unique word $\lfloor \bft\rfloor \in (\mathcal{X} \cup \mathcal{X}^*)^+$.
For instance, $\lfloor  x(x^2(yx^*)^*)^* zy^*\rfloor = xy(x^*)^3zy^*$.

\begin{remark}\label{rem:factor}
For any subterm $\bfs$ of a term $\bft$, either $\lfloor \bfs \rfloor$ or $\lfloor \bfs^* \rfloor$ is a factor of $\lfloor \bft\rfloor $.
\end{remark}

An \textit{identity} is an expression $\bfs \approx \bft$ formed by nonempty terms $\bfs, \bft \in \mathsf{T}(\mathcal{X})$, a \textit{word identity} is an identity $\bfu \approx \bfv$ formed by words $\bfu, \bfv \in (\mathcal{X} \cup \mathcal{X}^*)^+$.
We write $\bfu = \bfv$ if $\bfu$ and $\bfv$ are identical.  An identity $\bfu \approx \bfv$ is \textit{non-trivial} if $\bfu \neq \bfv$. An identity
$\bfs \approx \bft$ is directly deducible from an identity $\bfp \approx \bfq$ if there exists some substitution $\varphi : \mathcal{X}\rightarrow \mathsf{T}(\mathcal{X})$ such that $\varphi(\bfp)$ is a subterm of $\bfs$, and replacing this particular subterm
$\varphi(\bfp)$ of $\bfs$ with $\varphi(\bfq)$ results in the term $\bft$. An identity $\bfs \approx \bft$ is deducible from some set $\Sigma$
of identities if there exists a sequence $\bfs = \bfs_1,\bfs_2,\cdots,\bfs_r = \bft$ of terms such that each
identity $\bfs_i \approx \bfs_{i+1}$ is directly deducible from some identity in $\Sigma$.

\begin{remark}[{\cite[Sublemma~2.2]{CDE00}}]\label{rem:s=t}
 An identity $\bfs \approx \bft$ is deducible from \eqref{id: inv} if and only if $\lfloor \bfs \rfloor=\lfloor \bft \rfloor$.
\end{remark}

An  {\insem} $(S,\op)$ \textit{satisfies} an identity $\bfs \approx \bft$, if for any substitution $\varphi: \mathcal{X} \to S$, the elements $\varphi(\bfs)$ and $\varphi(\bft)$ of $S$ coincide; in this case, $\bfs \approx \bft$ is also said to be an \textit{identity of}  $(S,\op)$.
%Clearly any  {\inmon} that satisfies a  word identity $\bfs \approx \bft$ also satisfies the identity $\bfs [x_1, x_2, \ldots, x_n] \approx \bft[x_1, x_2, \ldots, x_n]$ for any distinct variables $\overline{x_1}, \overline{x_2}, \ldots, \overline{x_n} \in \mathcal{X}$.
%
%\begin{remark}\label{rem:delete}
%Let $\varphi$ be a substitution into an {\inmon} and $\bfu\in (\mathcal{X}\cup\mathcal{X}^*)^{+}$.
%Note that if $\varphi$ maps a variable $x$ in $\bfu$ to the unit element $1$ and otherwise to themselves, then $\varphi(\bfu)$ is the same as the image of the word obtained from $\bfu$ by removing all occurrences of $x$ and $x^*$ under $\varphi$.
%\end{remark}
\begin{remark}\label{rem:delete}
Note that assigning the unit element to a variable~$x$ in a word identity is effectively the same as removing all occurrences of~$x$ and~$x^*$. Therefore any {\inmon} that satisfies a word identity $\bfs \approx \bft$ also satisfies the word identity $\bfs [x_1, x_2, \ldots, x_n] \approx \bft[x_1, x_2, \ldots, x_n]$ for any distinct variables $\overline{x_1}, \overline{x_2}, \ldots, \overline{x_n} \in \mathcal{X}$.
\end{remark}

For any  {\insem} $(S,\op)$,
a set $\Sigma$ of identities of $(S,\op)$ is an \textit{identity basis} for  $(S,\op)$ if every identity satisfied by  $(S,\op)$ is deducible from $\Sigma$.
An  {\insem} is \textit{\fb} if it has some finite identity basis; otherwise, it is \textit{\nfb}.

The variety generated by semigroup $S$ [resp. {\insem} $(S,\op)$] is denoted by $\var S$ [resp. $\var (S,\op)$]. For any set $\Sigma$ of identities, denote by $\var \Sigma$ the variety determined by $\Sigma$.

\subsection{The Baxter monoid and its involution}%
Let $\mathcal{A} = \{ 1 < 2 < 3 < \cdots \}$ denote the set of positive integers, viewed as an infinite ordered alphabet.
The combinatorial objects and insertion algorithm related to the Baxter monoid are given in the following.

A \textit{right strict binary search tree} is a labelled rooted binary tree where the label of each
node is greater than or equal to the label of every node in its left subtree, and strictly
less than every node in its right subtree.
The associated insertion algorithm is as follows:
\begin{algo}\label{algo:sylv}
Input a right strict binary search tree $T$ and a symbol $a\in \mathcal{A}$.
If $T$ is empty, create a node and label it $a$. If $T$ is non-empty, examine the label $x$
of the root node: if $a > x$, recursively insert $a$ into the right subtree of the root node;
otherwise recursively insert $a$ into the left subtree of the root node. Output the resulting
tree.
\end{algo}

Let $w_1, \cdots, w_k\in \mathcal{A}$ and $\bfw=w_1\cdots w_k \in \mathcal{A}^{\star}$. Then the combinatorial object ${\rm P}_{\sylv_{\infty}}(\bfw)$ of $\bfw$ is obtained as follows: reading
$\bfw$ from right-to-left, one starts with an empty tree and inserts each symbol in $\bfw$ into a right strict binary search tree according to  Algorithm \ref{algo:sylv}.
For example, ${\rm P}_{\sylv_{\infty}}(36131512665)$ is given as follows:
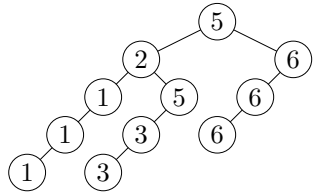
\begin{equation*}
\begin{tikzpicture}
[grow'=down,line width = 0pt,
every node/.style={draw,circle,inner sep=2pt},
level distance=.5cm,
level 1/.style={sibling distance=20mm},
level 2/.style={sibling distance=10mm},
level 3/.style={sibling distance=10mm},
level 4/.style={sibling distance=10mm}]

\node  (root) {5}
  child {node {6}
    child[missing]
    child {node {6}
       child[missing]
       child {node {6}}
      }
     }
  child {node {2}
    child {node {5}
      child[missing]
      child {node {3}
        child[missing]
        child {node {3}}
      }
    }
    child {node {1}
      child[missing]
      child {node{1}
        child[missing]
        child {node {1}}
      }
    }
  };
\end{tikzpicture}
\end{equation*}

A \textit{left strict binary search tree is} a labelled rooted binary tree where the label of each
node is strictly greater than the label of every node in its left subtree, and less than or
equal to every node in its right subtree.
The associated insertion algorithm is as follows:
\begin{algo}\label{algo:sylv-2}
Input a left strict binary search tree $T$ and a symbol $a \in \mathcal{A}$.
If $T$ is empty, create a node and label it $a$. If $T$ is non-empty, examine the label $x$
of the root node: if $a < x$, recursively insert $a$ into the left subtree of the root node;
otherwise recursively insert $a$ into the right subtree of the root note. Output the resulting
tree.
\end{algo}

Let $w_1, \cdots, w_k\in \mathcal{A}$ and $\bfw=w_1\cdots w_k \in \mathcal{A}^{\star}$. Then the combinatorial object ${\rm P}_{\sylv^{\sharp}_{\infty}}(\bfw)$ of $\bfw$ is obtained as follows: reading $\bfw$ from left-to-right, one starts with an empty tree and inserts each symbol in $\bfw$ into a left strict binary search tree according to  Algorithm \ref{algo:sylv-2}.
For example, ${\rm P}_{\sylv^{\sharp}_{\infty}}(36131512665)$ is given as follows:
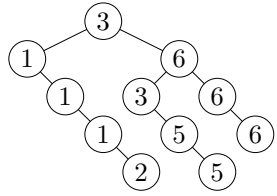
\begin{equation*}
\begin{tikzpicture}
[grow'=down,line width = 0pt,
every node/.style={draw,circle,inner sep=2pt},
level distance=.5cm,
level 1/.style={sibling distance=20mm},
level 2/.style={sibling distance=10mm},
level 3/.style={sibling distance=10mm},
level 4/.style={sibling distance=10mm}]

\node  (root) {3}
  child {node {6}
    child{node{6}
       child {node {6}}
       child [missing]}
    child {node {3}
      child {node {5}
        child {node {5}}
        child [missing]
        }
      child [missing]
    }
  }
  child {node {1}
    child {node {1}
      child {node {1}
        child {node {2}}
        child[missing]
        }
      child[missing]
      }
     child[missing]
     };
\end{tikzpicture}
\end{equation*}

Let $w_1, \cdots, w_k\in \mathcal{A}$ and $\bfw=w_1\cdots w_k \in \mathcal{A}^{\star}$. Then the combinatorial object  ${\rm P}_{\baxt_{\infty}}(\bfw)$ of $\bfw$ is obtained by Algorithms \ref{algo:sylv} and \ref{algo:sylv-2}, that is, ${\rm P}_{{\baxt}_{\infty}}(\bfw)= ({\rm P}_{{\sylv_{\infty}^{\sharp}}}(\bfw), {\rm P}_{{\sylv}_{\infty}}(\bfw))$.

Define the relation $\equiv_{\baxt_{\infty}}$ by
\[
\bfu \equiv_{\baxt_{\infty}} \bfv \quad  \text{if and only if} \quad {\rm P}_{\baxt_{\infty}}(\bfu) = {\rm P}_{\baxt_{\infty}}(\bfv)
\]
for any $\bfu,\bfv \in \mathcal{A}^{\star}$. The relation $\equiv_{\baxt_{\infty}}$ is a congruence on $\mathcal{A}^{\star}$. The Baxter monoid $\baxt_{\infty}$ is the factor monoid $ \mathcal{A}^{\star}/_{\equiv_{\baxt_{\infty}}}$. The rank-$n$ analogue $\baxt_{n}$ is the factor monoid $ \mathcal{A}^{\star}_n/_{\equiv_{\baxt_{\infty}}}$, where the relation $\equiv_{\baxt_{\infty}}$ is naturally restricted to $\mathcal{A}^{\star}_n\times\mathcal{A}^{\star}_n$ and $\mathcal{A}_n = \{1 < 2 < \cdots < n\}$ is the set of the first $n$ natural numbers viewed as a finite ordered alphabet.
It follows from the definition of $\equiv_{\baxt_{\infty}}$ that each element $[\bfu]_{\equiv_{\baxt_{\infty}}}$
of the factor monoid $\baxt_{\infty}$ can be identified with the
combinatorial object ${\rm P}_{\baxt_{\infty}}(\bfu)$.
Clearly $\baxt_1$ is a free monogenic monoid $\langle 1 \rangle=\{\varepsilon, 1, 1^2, 1^3,\ldots\}$ and so it is commutative.
Note that
\begin{align}\label{M-order}
\baxt_1 \subset \baxt_2 \subset \cdots \subset \baxt_i \subset\baxt_{i+1}\subset \cdots \subset \baxt_{\infty}.
\end{align}

For any word $\bfu\in \mathcal{A}^{\star}$, the \textit{length} $|\bfu|$ of $\bfu$  is the number of symbols occurring in $\bfu$, and $|\bfu|_a$ is the number of times the symbol $a$ appearing in $\bfu$; the \textit{evaluation} of $\bfu$, denoted by $\ev(\bfu)$, is the infinite
tuple of non-negative integers, indexed by $\mathcal{A}$, whose $a$-th element is $|\bfu|_a$, thus this tuple describes the number of each symbol in $\mathcal{A}$ that appears in $\bfu$. It is immediate from the definition of the Baxter monoid that if $\bfu \equiv_{\baxt_\infty} \bfv$, then $\ev(\bfu) = \ev(\bfv)$, and hence it makes sense to define the evaluation of each element of Baxter monoid to be the evaluation
of any word representing it. The \textit{support} of a word $\bfu\in \mathcal{A}^{\star}$, denoted by $\supp(\bfu)$, is the set of letters that occur in $\bfu$. Note that $\ev(\bfu) = \ev(\bfv)$ implies that $\supp(\bfu) = \supp(\bfv)$.

Let $\bfu \in\mathcal{A}^{\star}$ and $a, b \in \supp(\bfu)$ with $a < b$. We say that $\bfu$ has a
$b$-$a$ \textit{right precedence of index $r$} for some $r\geq 1$ if, when reading $\bfu$ from right to left, $b$ occurs $r$ times before the first
occurrence of $a$ and, for any $c \in \supp(\bfu)$ such that $a < c < b$, $c$ does not occur
before the first occurrence of $a$.
We say that $\bfu$ has a $a$-$b$ \textit{left precedence of index $\ell$} if, when reading $\bfu$ from
left to right, $a$ occurs $\ell$ times before the first occurrence of $b$ and, for any $c \in \supp(\bfu)$ such
that $a < c < b$, $c$ does not occur before the first occurrence of $b$.
Note that, for any given $a \in \supp(\bfu)$,
there is at most one $b \in \supp(\bfu)$ such that $\bfu$ has a $b$-$a$ right precedence of index
$r$; on the other hand, $\bfu$ can have several right precedences of the
form $b$-$x$ for a fixed $b$; and that, for any given $b \in \supp(\bfu)$,
there is at most one $a \in \supp(\bfu)$ such that $\bfu$ has a $a$-$b$ left precedence of index
$\ell$; on the other hand, $\bfu$ can have several left precedences of the
form $a$-$x$ for a fixed $a$.
Denote by
\[
\rpi(\bfu)=\{(b\textrm{-}a, r): \text{$\bfu$ has a $b\textrm{-}a$ right precedence of index $r$}\}
\]
and
\[
\lpi(\bfu)=\{(a\textrm{-}b, \ell): \text{$\bfu$ has a $a\textrm{-}b$ left precedence of index $\ell$}\}.
\]
For example, $\rpi(36131512665)=\{(2\textrm{-}1,1), (5\textrm{-}2,1), (5\textrm{-}3,2)\}$ and $\lpi(36131512665)=\{(1\textrm{-}2, 3), (3\textrm{-}5,2), (3\textrm{-}6,1)\}$.

\begin{proposition}[{\cite[Corollary~2.11]{CMR21b}}]\label{pro:inversion}
For any $\bfu, \bfv \in \mathcal{A}^{\star}$, $\bfu \equiv_{\baxt_\infty} \bfv$ if and only if $\ev(\bfu)=\ev(\bfv)$, $\rpi(\bfu)=\rpi(\bfv)$ and $\lpi(\bfu)=\lpi(\bfv)$.
\end{proposition}

Note that if the word $\bfu$ has a  $b$-$a$ right precedence or a $a$-$b$ left precedence for $a<b$, then each word in $[\bfu]_{\baxt_\infty}$ has a $b$-$a$ right precedence or a $a$-$b$ left precedence.

The Baxter monoid  can also be defined by the presentation $\left\langle \mathcal{A} \mid \mathcal{R}_{\baxt_\infty} \right\rangle$, where
\begin{align*}
\mathcal{R}_{\baxt_\infty} =& \left\{ (c\bfu ad\bfv b,  c\bfu da\bfv b): a \leq b < c \leq d \right\}\\
& \cup \left\{ (b\bfu da\bfv c, b\bfu ad\bfv c): a < b \leq c < d \right\}.
\end{align*}

For each $n \in \mathbb{N}$, a presentation for the Baxter monoid of rank $n$ can be obtained by restricting generators and relations of the above presentation to generators in $\mathcal{A}_n$. Note that these relations are length-preserving.

If the Baxter monoid $\baxt_n$ under the unary operation $^*$ is an involution monoid, then the relation $\equiv_{\baxt_{\infty}}$  must be compatible with the involution operation, that is, if $\bfu\equiv_{\baxt_{\infty}}\bfv$, then $\bfu^*\equiv_{\baxt_{\infty}}\bfv^*$.

Let $\mathcal{A}_n^{\sharp} := \{1^{\sharp} > 2^{\sharp} > \cdots > n^{\sharp}\}$ be the alphabet $\mathcal{A}_n$ on which the order relation has been reversed and ${(\mathcal{A}_n^{\sharp})}^{\sharp} := \mathcal{A}_n$. For $w_1w_2\cdots w_n \in \mathcal{A}^{\star}$, $(w_1w_2\cdots w_n)^{\sharp} := w_{n}^{\sharp} \cdots w_2^{\sharp} w_1^{\sharp}$.

\begin{proposition}[{\cite[Proposition~3.4]{Gir12}}]\label{pro:involution}
Let $\bfw$ and $\bfw'$ be two words in $\mathcal{A}_n^{\star}$. Then
$\bfw \equiv_{\baxt_{\infty}} \bfw'$ if and only if $\bfw^{\sharp} \equiv_{\baxt_{\infty}} {(\bfw')}^{\sharp}$.
\end{proposition}
The relation $\equiv_{\baxt_{\infty}}$ is compatible with the unary operation $^{\sharp}$. Thus  $(\baxt_n, ~^{\sharp})$ is an involution monoid, and
the involution $^{\sharp}$ is always called \textit{Sch\"{u}tzenberger's involution}.

\begin{proposition}
Sch\"{u}tzenberger's involution is the unique involution on the Baxter monoid $\baxt_n$ for each $n\geq 1$.
\end{proposition}
\begin{proof}
Suppose that $^*$ is an involution operation on $\baxt_n$. Note that the relations $\mathcal{R}_{\hypo_\infty}$ are length-preserving. Then the involution of a generator in $\mathcal{A}_n$ is still a generator in $\mathcal{A}_n$.
Since $\baxt_1$ has only one generator $1$, we have $1^*=1$. Thus the involution on $\baxt_1$ is trivial. For $\baxt_n$ with $n\geq 2$, let $a< b\leq n$. Then $ (b\bfu ab\bfv a)^* =a^*\bfv^*b^*a^*\bfu^*b^* \equiv_{\baxt_{\infty}}   a^*\bfv^*a^*b^*\bfu^*b^*=(b\bfu ba\bfv a)^*$ by $b\bfu ab\bfv a \equiv_{\baxt_{\infty}}  b\bfu ba\bfv a \in \mathcal{R}_{\hypo_\infty}$. This implies $ a^*\bfv^*b^*a^*\bfu^*b^* \equiv_{\baxt_{\infty}}   a^*\bfv^*a^*b^*\bfu^*b^* \in \mathcal{R}_{\baxt_\infty}$, and so $b^*< a^*$. Hence for any $a<b$, we must have $b^*< a^*$, whence $^*$ must be the order-reversing permutation on $\mathcal{A}_n$. Therefore $^*$ is Sch\"{u}tzenberger's involution, and so Sch\"{u}tzenberger's involution is the unique involution on the Baxter monoid $\baxt_n$.
\end{proof}

\subsection{Matrix representations over semirings}
Recall that $\mathbb{S}=(S, +, \cdot)$ is a \textit{commutative semiring} with additive identity $\mathbf{0}$ and multiplicative identity $\mathbf{1}$ if $S$ is a set equipped with two binary operations $+$ and $\cdot$ such that $(S, +)$ and $(S, ~\cdot)$ are  commutative monoids satisfying
\[
a(b+c)=a\cdot b+a\cdot c \quad\text{and}\quad \mathbf{0}\cdot a=\mathbf{0}
\]
for all $a,b,c \in \mathbb{S}$. Semiring $\mathbb{S}$ is \textit{idempotent} if $a+a=a$ for all $a\in \mathbb{S}$.
An element $a\in \mathbb{S}$ is \textit{infinite multiplicative order} if for any non-negative integers $i,j$, $a^{i}=a^{j}$ if and only if $i=j$.
In this paper, we always assume that $\mathbb{S}$ is a commutative and idempotent semiring with $\mathbf{0}, \mathbf{1}$ containing an element of infinite multiplicative order.
A common example of such semiring is the tropical semiring
$\mathbb{T}= (\mathbb{R} \cup \{-\infty\}, \oplus, \otimes)$, which is the set $\mathbb{R}$ of real numbers together with minus infinity $-\infty$, with the addition and multiplication defined as follows
\[
a \oplus b = \max \{a, b\} \quad \mbox{and} \quad a \otimes b = a + b.
\]
In other words, the tropical sum of two numbers is their maximum and the tropical
product of two numbers is their sum, and
$-\infty$ is the additive identity and $0$ is the multiplicative identity. Note that except $-\infty$ and $0$, all other elements in $\mathbb{T}$ have infinite multiplicative order.

Note that the set of all $n \times n$  matrices with entries in $\mathbb{S}$ forms a semigroup under the matrix multiplication induced from the operations in $\mathbb{S}$. We denote this semigroup by $M_{n}(\mathbb{S})$. Let $UT_n(\mathbb{S})$ be the subsemigroup of $M_n(\mathbb{S})$ of all upper triangular $n \times n$  matrices.
For any matrix $A\in M_n(\mathbb{S})$, denote by $A^{D}$ the matrix obtained by reflecting $A$ with respect to the secondary diagonal (from the top right to the bottom left corner), that is,  $(A^{D})_{ij}=A_{(n+1-j)(n+1-i)}$. It is easy to verify that this unary operation $^D$ is an involution operation of $M_{n}(\mathbb{S})$.
A (linear) representation of a semigroup $S$ [resp. involution semigroup $(S,~^*)$] is a
homomorphism $\rho : S \rightarrow M_n(\mathbb{S})$ [resp. $\rho : (S,~^*) \rightarrow (M_n(\mathbb{S}),~^{D})$ ]. The homomorphism $\rho$ is said to be \textit{faithful} if it is injective. Note that an involution semigroup representation $\rho : (S,~^*) \rightarrow (M_n(\mathbb{S}),~^{D})$ deduces a semigroup representation $\rho : S \rightarrow M_n(\mathbb{S})$, but a semigroup representation $\rho : S \rightarrow M_n(\mathbb{S})$ can not be an involution semigroup representation $\rho : (S,~^*) \rightarrow (M_n(\mathbb{S}),~^{D})$. The tropical semiring is of interest as a natural carrier for representations of semigroups. For example,  the bicyclic monoid $\mathcal{B}:=\langle a,b\mid ba=1\rangle$, which is ubiquitous in infinite semigroup theory, admits no faithful finite dimensional representations over any field; however it has a number of natural representations over the tropical semiring \cite{DJK18,Izhakian10}.

\section{matrix representations of $(\mathsf{baxt}_n, ~^\sharp)$}\label{sec:repre}%
In this section, we exhibit a faithful representation of $(\mathsf{baxt}_n, ~^\sharp)$ for each finite $n$ as an involution monoid of upper triangular matrices over $\mathbb{S}$ under the skew transposition, and we prove that all involution semigroups $(\baxt_n, ~^{\sharp})$ with $n\geq 4$ generate the same variety.

For convenience, denote by
 \begin{gather*}
 \mathrm{P}=\begin{pmatrix}
s &  \mathbf{0} \\
\mathbf{0} &  \mathbf{1}
\end{pmatrix},
\mathrm{Q}=\begin{pmatrix}
\mathbf{1} & \mathbf{0} \\
\mathbf{0} & s
\end{pmatrix},
 \mathrm{J}=\begin{pmatrix}
\mathbf{1} & \mathbf{1}  \\
\mathbf{0} & \mathbf{0}
\end{pmatrix},
\mathrm{K}=\begin{pmatrix}
\mathbf{0} &\mathbf{1}  \\
\mathbf{0} & \mathbf{1}
\end{pmatrix}
\end{gather*}
where $s\in \mathbb{S}$ is an element of infinite multiplicative order. Denote by $\mathsf{diag}\{\Lambda_1,\Lambda_2,\dots, \Lambda_n\}$ the block diagonal matrix
\[
\begin{pmatrix}
\Lambda_1 & & & \\
 & \Lambda_2 & &\\
 & & \ddots &\\
  &&& \Lambda_n
\end{pmatrix}
\]
where $\Lambda_1,\Lambda_2,\dots, \Lambda_n$ are square matrices. And let $\mathrm{E}_n$ be the $n\times n$ matrix with $\mathbf{1}$s on the main diagonal and $\mathbf{0}$s elsewhere.

First we give a matrix representation of $(\mathsf{baxt}_1,~^\sharp)$.
Define a map $\varphi_1: \mathcal{A}_1\cup \{\varepsilon\}\rightarrow UT_2(\mathbb{S})$ given by $\varepsilon\mapsto \mathrm{E}_2$ and
$1\mapsto \mathrm{P}\mathrm{Q}$.
Clearly, the map $\varphi_1$ induces a faithful representation of $(\mathsf{baxt}_1,~^\sharp)$.

Next we consider the matrix representation of $(\mathsf{baxt}_2,~^\sharp)$.
Define a map $\varphi_2: \mathcal{A}_2\cup \{\varepsilon\}\rightarrow UT_6(\mathbb{S})$ given by $\varepsilon\mapsto \mathrm{E}_6,$
\[
1\mapsto \mathsf{diag}\{s,\mathrm{P}, \mathrm{J}, \mathbf{1}\},\;\; 2\mapsto\mathsf{diag}\{\mathbf{1},\mathrm{K}, \mathrm{Q}, s\}.
\]
Clearly, $\varphi_2$ can be extended to a homomorphism from $\mathcal{A}_2^{\star}$ to $UT_6(\mathbb{S})$. Note that
$\varphi_2(1^{\sharp})= \varphi_2(2)=\mathsf{diag}\{\mathbf{1},\mathrm{K}, \mathrm{Q}, s\}=(\varphi_2(1))^{D}$ and $\varphi_2(2^{\sharp})= \varphi_2(1)=\mathsf{diag}\{s,\mathrm{P}, \mathrm{J}, \mathbf{1}\}=(\varphi_2(2))^{D}$.
Thus $\varphi_2$ can be extended to a homomorphism from $(\mathcal{A}_2^{\star}, ~^{\sharp})$ to $(UT_6(\mathbb{S}),~^{D})$.
In fact, $\varphi_2$ induces a faithful representation of $(\mathsf{baxt}_2,~^\sharp)$.

\begin{theorem}\label{thm:baxt2-repre}
The map $\varphi_2:(\mathsf{baxt}_2,~^\sharp)\rightarrow(UT_6(\mathbb{S}),~^{D})$ is a faithful representation of $(\mathsf{baxt}_2,~^\sharp)$.
\end{theorem}

\begin{proof}
Note that $\varphi_2$ is a homomorphism from $(\mathcal{A}_2^{\star}, ~^{\sharp})$ to $(UT_6(\mathbb{S}),~^{D})$. Then to show that the map $\varphi_2$ induces a homomorphism from $(\mathsf{baxt}_2,~^\sharp)$ to $(UT_6(\mathbb{S}),~^{D})$, we only need to show that for any $\bfu, \bfv\in \mathcal{A}_2^{\star}$, if $\bfu\equiv_{\baxt_\infty}\bfv$, then
$\varphi_2(\bfu) = \varphi_2(\bfv)$.
By the definition of $\varphi_2$, it is easy to verify that for any $\bfw\in \mathcal{A}_2^{\star}$,
\[
\varphi_2(\bfw)=\mathsf{diag}\{\Lambda_1,\Lambda_2,\Lambda_3, \Lambda_4\}
\]
where
\begin{align*}\label{id:2}
\Lambda_1&=
\left\{
  \begin{array}{ll}
s^{|\bfw|_1}, &\, \,\hbox{if $\{1\}\in \supp(\bfw)$,} \\ [0.1cm]
\mathbf{1}, & \,\,\hbox{if $\{1\}\not\in \supp(\bfw)$,}
\end{array}
\right.\\
\Lambda_2&=
\left\{
  \begin{array}{ll}
\mathrm{P}^{|\bfw|_1}, & \hbox{if $\supp(\bfw)=\{1\}$,} \\ [0.1cm]
\mathrm{K}, & \hbox{if $\{2\}\subseteq \supp(\bfw)$ and $(1\textrm{-}2, \ell) \not\in \lpi(\bfw)$ for any $\ell> 0$,} \\ [0.1cm]
\mathrm{P}^{\ell_1}\mathrm{K}, & \hbox{if $\supp(\bfw)=\{1,2\}$ and $(1\textrm{-}2, \ell_1)\in \lpi(\bfw)$,}
\end{array}
\right.\\
\Lambda_3&=
\left\{
  \begin{array}{ll}
\mathrm{Q}^{|\bfw|_2}, & \hbox{if $\supp(\bfw)=\{2\}$,} \\ [0.1cm]
\mathrm{J}, & \hbox{if $\{1\}\subseteq \supp(\bfw)$ and $(2\textrm{-}1, r) \not\in \rpi(\bfw)$ for any $r> 0$,} \\ [0.1cm]
\mathrm{J}\mathrm{Q}^{r_1}, & \hbox{if $\supp(\bfw)=\{1,2\}$ and $(2\textrm{-}1, r_1)\in \rpi(\bfw)$,}
\end{array}
\right.\\
\Lambda_4&=
\left\{
  \begin{array}{ll}
s^{|\bfw|_2}, &\,\, \hbox{if $\{2\}\in \supp(\bfw)$,} \\ [0.1cm]
\mathbf{1}, &\,\, \hbox{if $\{2\}\not\in \supp(\bfw)$.}
\end{array}
\right.
\end{align*}
Since $\ev(\bfu)=\ev(\bfv), \lpi(\bfu)=\lpi(\bfv), \rpi(\bfu)=\rpi(\bfv)$ by Proposition \ref{pro:inversion}, it is routine to show that $\varphi_2(\bfu) = \varphi_2(\bfv)$.

Suppose that $\bfu\not\equiv_{\baxt_\infty}\bfv$. Then $\ev(\bfu)\neq\ev(\bfv)$, or $\lpi(\bfu)\neq\lpi(\bfv)$, or $\rpi(\bfu)\neq\rpi(\bfv)$ by Proposition \ref{pro:inversion}. By the definition of $\varphi_2$, it is routine to show that $\varphi_2(\bfu)\neq\varphi_2(\bfv)$.
Hence $\varphi_2$ is injective. Therefore the map $\varphi_2:(\mathsf{baxt}_2, ~^\sharp)\rightarrow(UT_6(\mathbb{S}),~^{D})$ is a faithful representation of $(\mathsf{baxt}_2,~^\sharp)$.
\end{proof}

Next we consider the matrix representation of $(\mathsf{baxt}_3,~^\sharp)$.
Define a map $\varphi_3: \mathcal{A}_3\cup \{\varepsilon\}\rightarrow UT_{15}(\mathbb{S})$ given by
$\varepsilon\mapsto \mathrm{E}_{15}$,
\begin{align*}
1&\mapsto
\mathsf{diag}\{s, \mathrm{P},\mathrm{P},\mathrm{E}_2, \mathbf{1}, \mathrm{J},\mathrm{E}_2, \mathrm{J}, \mathbf{1}\},\\
2&\mapsto
\mathsf{diag}\{\mathbf{1}, \mathrm{K},\mathrm{K},\mathrm{P},s, \mathrm{Q},\mathrm{J}, \mathrm{J},\mathbf{1}\},\\
3&\mapsto
\mathsf{diag}\{\mathbf{1}, \mathrm{K},\mathrm{E}_2, \mathrm{K},\mathbf{1}, \mathrm{E}_2,\mathrm{Q},\mathrm{Q}, s\}.
\end{align*}
Clearly, $\varphi_3$ can be extended to a homomorphism from $\mathcal{A}_3^{\star}$ to $UT_{15}(\mathbb{S})$. Note that $\varphi_3(1^{\sharp})=\varphi_3(3)=(\varphi_3(1))^{D}, \varphi_3(2^{\sharp})=\varphi_3(2)=(\varphi_3(2))^{D}$ and
$\varphi_3(3^{\sharp})=\varphi_3(1)=(\varphi_3(1))^{D}$.
Thus $\varphi_3$ can be extended to a homomorphism from $(\mathcal{A}_3^{\star}, ~^{\sharp})$ to $(UT_{15}(\mathbb{S}),~^{D})$.
In fact, $\varphi_3$ induces a faithful representation of  $(\mathsf{baxt}_3,~^\sharp)$.

\begin{theorem}\label{thm:baxt3-repre}
The map $\varphi_3:(\mathsf{baxt}_3,~^\sharp)\rightarrow(UT_{15}(\mathbb{S}),~^{D})$ is a faithful representation of $(\mathsf{baxt}_3,~^\sharp)$.
\end{theorem}

\begin{proof}
Note that $\varphi_3$ is a homomorphism from $(\mathcal{A}_3^{\star}, ~^{\sharp})$ to $(UT_{15}(\mathbb{S}),~^{D})$. Then to show that the map $\varphi_3$ induces a homomorphism from $(\mathsf{baxt}_3,~^\sharp)$ to $(UT_{15}(\mathbb{S}),~^{D})$, we only need to show that for any $\bfu,\bfv\in \mathcal{A}_3^{\star}$, if $\bfu\equiv_{\hypo_\infty}\bfv$, then $\varphi_3(\bfu) = \varphi_3(\bfv)$.
By the definition of $\varphi_3$, it is easy to verify that for any $\bfw\in \mathcal{A}_3^{\star}$,
\[
\varphi_3(\bfw)=\mathsf{diag}\{\Lambda_1,\Lambda_2,\Lambda_3,\Lambda_4,\Lambda_5, \Lambda_6,\Lambda_7, \Lambda_8, \Lambda_9\}
\]
where
\begin{align*}\label{id:3}
\Lambda_1&=
\left\{
  \begin{array}{ll}
s^{|\bfw|_1}, & \,\,\hbox{if $\{1\}\in \supp(\bfw)$,} \\ [0.1cm]
\mathbf{1}, & \,\,\hbox{if $\{1\}\not\in \supp(\bfw)$,}
\end{array}
\right.\\
\Lambda_2&= \left\{
  \begin{array}{ll}
\mathrm{P}^{|\bfw|_1}, & \hbox{if $\supp(\bfw)=\{1\}$,} \\ [0.1cm]
%\mathrm{K}, & \hbox{if $\supp(\bfw)=\{2\}$, $\{3\}$  or $\{2, 3\}$,} \\ [0.1cm]
\mathrm{P}^{\ell_1}\mathrm{K}, & \hbox{if $\{1, 2\}= \supp(\bfw)$ and $(1\textrm{-}2, \ell_1) \in \lpi(\bfw)$,} \\ [0.1cm]
%\mathrm{K}, & \hbox{if $\{1, 2\}\subseteq \supp(\bfw)$ and $(1\textrm{-}2, k) \not\in \lpi(\bfw)$ for any $k>0$,} \\ [0.1cm]
\mathrm{P}^{\ell_2}\mathrm{K}, & \hbox{if $\{1, 3\}\subseteq \supp(\bfw)$ and $(1\textrm{-}3, \ell_2) \in \lpi(\bfw)$,}\\ [0.1cm]
%\mathrm{K}, & \hbox{if $\{1, 3\}= \supp(\bfw)$ and $(1\textrm{-}3, k) \not\in \lpi(\bfw)$ for any $k> 0$,}\\ [0.1cm]
\mathrm{P}^{\ell_1}\mathrm{K}, & \hbox{if $\{1, 2, 3\}= \supp(\bfw)$ and $(1\textrm{-}2, \ell_1)\in \lpi(\bfw)$, $(2\textrm{-}3, \ell_3)\in \lpi(\bfw)$,}\\ [0.1cm]
\mathrm{K}, & \hbox{otherwise,}
\end{array}
\right.\\
\Lambda_3&= \left\{
\begin{array}{ll}
\mathrm{P}^{|\bfw|_1}, & \hbox{if $\supp(\bfw)=\{1\}$ or $\{1,3 \}$,} \\ [0.1cm]
%\mathrm{K}, & \hbox{if $\supp(\bfw)=\{2\}$ or $\{2, 3\}$,} \\ [0.1cm]
\mathrm{E}_2, & \hbox{if $\supp(\bfw)=\{3\}$,} \\ [0.1cm]
\mathrm{P}^{\ell_1}\mathrm{K}, & \hbox{if $\{1, 2\}\subseteq \supp(\bfw)$ and $(1\textrm{-}2, \ell_1) \in \lpi(\bfw)$,} \\ [0.1cm]
\mathrm{K}, & \hbox{otherwise,}
\end{array}
\right.\\
\Lambda_4&= \left\{
  \begin{array}{ll}
\mathrm{E}_2, & \hbox{if $\supp(\bfw)=\{1\}$,} \\ [0.1cm]
\mathrm{P}^{|\bfw|_2}, & \hbox{if $\supp(\bfw)=\{2\}$ or $\{1,2\}$,} \\ [0.1cm]
%\mathrm{K}, & \hbox{if $\supp(\bfw)=\{3\}$ or $\{1,3\}$,} \\ [0.1cm]
\mathrm{P}^{\ell_3}\mathrm{K},  & \hbox{if $\{2,3\}\subseteq\supp(\bfw)$ and $(2\textrm{-}3, \ell_3) \in \lpi(\bfw)$,}\\ [0.1cm]
\mathrm{K},  & \hbox{otherwise,}
\end{array}
\right.\\
\Lambda_5& =
\left\{
  \begin{array}{ll}
s^{|\bfw|_2}, &\,\, \hbox{if $\{2\}\in \supp(\bfw)$,} \\ [0.1cm]
\mathbf{1}, &\,\, \hbox{if $\{2\}\not\in \supp(\bfw)$,}
\end{array}
\right.\\
\Lambda_6 &= \left\{
  \begin{array}{ll}
%\mathrm{J}, & \hbox{if $\supp(\bfw)=\{1\}$ or $\{1, 3\}$,} \\ [0.1cm]
\mathrm{Q}^{|\bfw|_2}, & \hbox{if $\supp(\bfw)=\{2\}$ or $\{2, 3\}$,} \\ [0.1cm]
\mathrm{E}_2, & \hbox{if $\supp(\bfw)=\{3\}$,} \\ [0.1cm]
\mathrm{J}\mathrm{Q}^{r_1}, & \hbox{if $\{1, 2\}\subseteq \supp(\bfw)$ and $(2\textrm{-}1, r_1) \in \rpi(\bfw)$,} \\ [0.1cm]
\mathrm{J}, & \hbox{otherwise,}
\end{array}
\right.\\
\Lambda_7 &= \left\{
\begin{array}{ll}
\mathrm{E}_2, & \hbox{if $\supp(\bfw)=\{1\}$,} \\ [0.1cm]
%\mathrm{J}, & \hbox{if $\supp(\bfw)=\{2\}$ or $\{1, 2\}$,} \\ [0.1cm]
\mathrm{Q}^{|\bfw|_3}, & \hbox{if $\supp(\bfw)=\{3\}$ or $\{1, 3\}$,} \\ [0.1cm]
\mathrm{J}\mathrm{Q}^{r_3}, & \hbox{if $\{2, 3\}\subseteq \supp(\bfw)$ and $(3\textrm{-}2, r_3) \in \rpi(\bfw)$,} \\ [0.1cm]
\mathrm{J}, & \hbox{otherwise,}
\end{array}
\right.\\
\Lambda_8 &= \left\{
  \begin{array}{ll}
%\mathrm{J}, & \hbox{if $\supp(\bfw)=\{1\}$, $\{2 \}$ or $\{1, 2\}$,} \\ [0.1cm]
\mathrm{Q}^{|\bfw|_3}, & \hbox{if $\supp(\bfw)=\{3\}$,} \\ [0.1cm]
\mathrm{J}\mathrm{Q}^{r_2},  & \hbox{if $\{1,3\}\subseteq\supp(\bfw)$ and $(3\textrm{-}1, r_2) \in \rpi(\bfw)$,}\\ [0.1cm]
%\mathrm{J},  & \hbox{if $\{1,3\}\subseteq\supp(\bfw)$ and $(3\textrm{-}1, r) \not\in \rpi(\bfw)$ for any $r\geq 0$,}\\ [0.1cm]
\mathrm{J}\mathrm{Q}^{r_3},  & \hbox{if $\{2,3\}=\supp(\bfw)$ and $(3\textrm{-}2, r_3) \in \rpi(\bfw)$,}\\ [0.1cm]
\mathrm{J}\mathrm{Q}^{r_3},  & \hbox{if $\{1, 2,3\}=\supp(\bfw)$ and $(2\textrm{-}1, r_1) \in \rpi(\bfw)$, $(3\textrm{-}2, r_3) \in \rpi(\bfw)$,}\\ [0.1cm]
\mathrm{J},  & \hbox{otherwise,}
\end{array}
\right.\\
\Lambda_9&=
\left\{
  \begin{array}{ll}
s^{|\bfw|_3}, &\,\, \hbox{if $\{3\}\in \supp(\bfw)$,} \\ [0.1cm]
\mathbf{1}, &\,\, \hbox{if $\{3\}\not\in \supp(\bfw)$.}
\end{array}
\right.
\end{align*}
Since $\ev(\bfu)=\ev(\bfv), \lpi(\bfu)=\lpi(\bfv), \rpi(\bfu)=\rpi(\bfv)$ by Proposition \ref{pro:inversion}, it is routine to show that $\varphi_3(\bfu) = \varphi_3(\bfv)$.

Suppose that $\bfu\not\equiv_{\baxt_\infty}\bfv$. Then $\ev(\bfu)\neq\ev(\bfv)$, or $\rpi(\bfu)\neq\rpi(\bfv)$, or $\rpi(\bfu)\neq\rpi(\bfv)$ by Proposition \ref{pro:inversion}. By the definition of $\varphi_3$, it is routine to show that $\varphi_3(\bfu)\neq\varphi_3(\bfv)$.
Hence $\varphi_3$ is injective. Therefore the map $\varphi_3:(\mathsf{baxt}_3, ~^\sharp)\rightarrow(UT_{15}(\mathbb{S}),~^{D})$ is a faithful representation of $(\mathsf{baxt}_3,~^\sharp)$.
\end{proof}

Now we consider the matrix representation of $(\mathsf{baxt}_n,~^\sharp)$ for $n\geq 4$.
For any $([\bfu]_{\baxt_3}, [\bfv]_{\baxt_3})\in \baxt_3\times \baxt_3$, define an involution operation $^{\sharp}$ on $\baxt_3$ $\times \baxt_3$ by  $([\bfu]_{\baxt_3}, [\bfv]_{\baxt_3})^{\sharp}= ([\bfv]_{\baxt_3}^{\sharp}, [\bfu]_{\baxt_3}^{\sharp})$. For any $i<j \in \mathcal{A}_n$ with $n\geq 4$, by the definition of $^{\sharp}$, we have that $j^{\sharp}<i^{\sharp}$ and there is at most one $k\in\mathcal{A}_n$ satisfying $k=k^{\sharp}$ and $j-i= i^{\sharp}-j^{\sharp}$ when $i\neq i^{\sharp}, j\neq j^{\sharp}$. So there are seven cases about the order of $i,j, i^{\sharp}, j^{\sharp}$ in $\mathcal{A}_n$: $i^{\sharp}=j$,  $i<j=j^{\sharp}<i^{\sharp}$, $j^{\sharp}<i=i^{\sharp}<j$, $i<j<j^{\sharp}<i^{\sharp}$, $j^{\sharp}<i^{\sharp}<i<j$, $i<j^{\sharp}<j<i^{\sharp}$, $j^{\sharp}<i<i^{\sharp}<j$.
For any $i<j \in \mathcal{A}_n$ with $n\geq 4$, define a map $\varphi_{ij}$ from $\mathcal{A}_n^{*}$ to $\baxt_3\times \baxt_3$ which can be determined by the following four cases according to the order of $i, i^{\sharp}, j, j^{\sharp}$ in $\mathcal{A}_n$.

{\bf Case 1.} $i^{\sharp} = j$. Define a map $\lambda: \mathcal{A}_n \rightarrow \baxt_3$ given by
\begin{align*}
k &\mapsto \begin{cases}
[1]_{\baxt_3} & \text{if}\ k = i,\\
[3]_{\baxt_3} & \text{if}\ k = j,\\
[31]_{\baxt_3} & \text{if}\ i < k < j,\\
\left[\varepsilon\right]_{\baxt_3} & \text{otherwise}.
\end{cases}
\end{align*}
Clearly, this map can be extended to a homomorphism $\lambda: \mathcal{A}_n^{\star} \rightarrow \baxt_3$.
Define a map $\lambda_{ij}: \mathcal{A}_n \rightarrow \baxt_3\times \baxt_3$ given by
\[
k \mapsto (\lambda(k), \lambda(k)).
\]
This map can be extended to a homomorphism $\lambda_{ij}: \mathcal{A}_n^{\star} \rightarrow \baxt_3\times \baxt_3$. Further, $\lambda_{ij}$
is also a homomorphism from $(\mathcal{A}_n^{\star},~^{\sharp})$ to $(\baxt_3\times \baxt_3,~^{\sharp})$. This is because for any $k \in \mathcal{A}_n$, $\lambda(k^{\sharp})=(\lambda(k))^{\sharp}$ which follows from
\begin{align*}
\begin{cases}
\lambda(k^{\sharp})=(\lambda(k))^{\sharp}=[3]_{\baxt_3} & \text{if}\ k = i,\\
\lambda(k^{\sharp})=(\lambda(k))^{\sharp}=[31]_{\baxt_3} & \text{if}\ i < k < j,\\
\lambda(k^{\sharp})=(\lambda(k))^{\sharp}=[1]_{\baxt_3}  & \text{if}\ k = j,\\
\lambda(k^{\sharp})=(\lambda(k))^{\sharp}=[\varepsilon]_{\baxt_3}   & \text{otherwise}.
\end{cases}
\end{align*}
Therefore for any $\bfw=k_1k_2\cdots k_n$,
\begin{align*}
\lambda_{ij}({\bfw}^{\sharp})&= (\lambda({\bfw}^{\sharp}), \lambda({\bfw}^{\sharp}))\\
&=(\lambda(k_n^{\sharp})\cdots\lambda(k_1^{\sharp}), \lambda(k_n^{\sharp})\cdots\lambda(k_1^{\sharp}))\\
&=((\lambda(k_n))^{\sharp}\cdots(\lambda(k_1))^{\sharp}, (\lambda(k_n))^{\sharp}\cdots(\lambda(k_1))^{\sharp})\\
&=((\lambda(\bfw))^{\sharp}, (\lambda(\bfw))^{\sharp})\\
&=(\lambda_{ij}(\bfw))^{\sharp}.
\end{align*}
%Note that $\lambda_{ij} (w)$ is the baxtplactic class of the word obtained from $w$ by replacing any occurrence of $i$ by $1$, any occurrence of $j$ by $3$, any occurrence of $k$ with $i < k < j$ by $31$, and erasing any occurrence of any other element.

{\bf Case 2.} $i<j=j^{\sharp}<i^{\sharp}$ or $j^{\sharp}<i=i^{\sharp}<j$. For convenience, let $i_1=i, i_2=j$ and $i_3=i^{\sharp}$ when $i<j=j^{\sharp}<i^{\sharp}$ and $i_1=j^{\sharp}, i_2=i$ and $i_3=j$ when $j^{\sharp}<i=i^{\sharp}<j$. Define maps $\theta_1: \mathcal{A}_n \rightarrow \baxt_3$ and $\theta_2: \mathcal{A}_n \rightarrow \baxt_3$ by
\begin{align*}
k \mapsto \begin{cases}
[1]_{\baxt_3} & \text{if}\ k = i_1,\\
[2]_{\baxt_3} & \text{if}\ k = i_2,\\
[21]_{\baxt_3} & \text{if}\ i_1 < k < i_2,\\
\left[\varepsilon\right]_{\baxt_3} & \text{otherwise,}
\end{cases}\;\;\; \mbox{and}\;\;\;
k \mapsto \begin{cases}
[2]_{\baxt_3} & \text{if}\ k = i_2,\\
[3]_{\baxt_3} & \text{if}\ k = i_3,\\
[32]_{\baxt_3} & \text{if}\ i_2 < k < i_3,\\
\left[\varepsilon\right]_{\baxt_3} & \text{otherwise,}
\end{cases}
\end{align*}
respectively.
Clearly, $\theta_1, \theta_2$ can be extended to homomorphisms from $\mathcal{A}_n^{\star}$ to $\baxt_3$ respectively.
Define a map $\theta_{ij}: \mathcal{A}_n \rightarrow \baxt_3\times \baxt_3$ by
\[
k \mapsto (\theta_1(k), \theta_2(k)).
\]
This map can be extended to a homomorphism $\theta_{ij}: \mathcal{A}_n^{\star} \rightarrow \baxt_3\times \baxt_3$. Further, $\theta_{ij}$ is also a homomorphism from $(\mathcal{A}_n^{\star},~^{\sharp})$ to $(\baxt_3\times \baxt_3,~^{\sharp})$. This is because for any $k \in \mathcal{A}_n$, $\theta_1(k^{\sharp})=(\theta_2(k))^{\sharp}, (\theta_1(k))^{\sharp}=\theta_2(k^{\sharp})$ which follows from
\begin{align*}
\begin{cases}
\theta_1(k^{\sharp})=(\theta_2(k))^{\sharp}=[\varepsilon]_{\baxt_3}, \;\;\;(\theta_1(k))^{\sharp}= \theta_2(k^{\sharp})=[3]_{\baxt_3} & \text{if}\ k = i_1,\\
\theta_1(k^{\sharp})=(\theta_2(k))^{\sharp}=[\varepsilon]_{\baxt_3}, \;\;\;(\theta_1(k))^{\sharp}= \theta_2(k^{\sharp})=[32]_{\baxt_3} & \text{if}\ i_1 < k < i_2,\\
\theta_1(k^{\sharp})=(\theta_2(k))^{\sharp}=[2]_{\baxt_3}, \;\;\;(\theta_1(k))^{\sharp}= \theta_2(k^{\sharp})=[2]_{\baxt_3}  & \text{if}\ k = i_2,\\
\theta_1(k^{\sharp})=(\theta_2(k))^{\sharp}=[21]_{\baxt_3},\;(\theta_1(k))^{\sharp}= \theta_2(k^{\sharp})=[\varepsilon]_{\baxt_3}& \text{if}\  i_2<k<i_3,\\
\theta_1(k^{\sharp})=(\theta_2(k))^{\sharp}=[1]_{\baxt_3}, \;\;\; (\theta_1(k))^{\sharp}=\theta_2(k^{\sharp})=[\varepsilon]_{\baxt_3} & \text{if}\ k =i_3,\\
\theta_1(k^{\sharp})=(\theta_2(k))^{\sharp}=[\varepsilon]_{\baxt_3}, \;\;\;(\theta_1(k))^{\sharp}= \theta_2(k^{\sharp})=[\varepsilon]_{\baxt_3} & \text{if $k<i_1$ or $k>i_3$}.
\end{cases}
\end{align*}
Therefore for any $\bfw=k_1k_2\cdots k_n$,
\begin{align*}
\theta_{ij}({\bfw}^{\sharp})&= (\theta_1({\bfw}^{\sharp}), \theta_2({\bfw}^{\sharp}))\\
&=(\theta_1(k_n^{\sharp})\cdots\theta_1(k_1^{\sharp}), \theta_2(k_n^{\sharp})\cdots\theta_2(k_1^{\sharp}))\\
&=((\theta_2(k_n))^{\sharp}\cdots(\theta_2(k_1))^{\sharp}, (\theta_1(k_n))^{\sharp}\cdots(\theta_1(k_1))^{\sharp})\\
&=((\theta_2(\bfw))^{\sharp}, (\theta_1(\bfw))^{\sharp})\\
&=(\theta_{ij}(\bfw))^{\sharp}.
\end{align*}
%Note that $\theta_{ij}(w)$ is a pair of baxtplactic class of the word obtained from $w$. The first component can be obtained from $w$ by replacing any occurrence of $i$ by $1$; any occurrence of $j$ by $2$; any occurrence of $k$ with $i < k < j$ by $21$; and erasing any occurrence of any other element. The second component can be obtained from $w$ by replacing any occurrence of $j$ by $2$; any occurrence of $i^{\sharp}$ by $3$; any occurrence of $k$ with $j < k < i^{\sharp}$ by $32$; and erasing any occurrence of any other element.

{\bf Case 3.} $i<j<j^{\sharp}<i^{\sharp}$ or $j^{\sharp}<i^{\sharp}<i<j$. For convenience, let $i_1=i, i_2=j, i_3=j^{\sharp}$ and $i_4=i^{\sharp}$ when $i<j<j^{\sharp}<i^{\sharp}$ and $i_1=j^{\sharp}, i_2=i^{\sharp},i_3=i$ and $i_4=j$ when $j^{\sharp}<i^{\sharp}<i<j$. Define maps $\eta_1: \mathcal{A}_n \rightarrow \baxt_3$ and $\eta_2: \mathcal{A}_n \rightarrow \baxt_3$ by
\begin{align*}
k \mapsto \begin{cases}
[1]_{\baxt_3} & \text{if}\ k = i_1,\\
[2]_{\baxt_3} & \text{if}\ k = i_2,\\
[21]_{\baxt_3} & \text{if}\ i_1 < k < i_2,\\
\left[\varepsilon\right]_{\baxt_3} & \text{otherwise},
\end{cases}\;\;\; \mbox{and}\;\;\;
k \mapsto \begin{cases}
[2]_{\baxt_3} & \text{if}\  k =i_3,\\
[3]_{\baxt_3} & \text{if}\  k = i_4,\\
[32]_{\baxt_3} & \text{if}\ i_3 < k < i_4,\\
\left[\varepsilon\right]_{\baxt_3} & \text{otherwise,}
\end{cases}
\end{align*}
respectively. Clearly, $\eta_1, \eta_2$ can be extended to homomorphisms from $\mathcal{A}_n^{\star}$ to $\baxt_3$ respectively.
Define a map $\eta_{ij}: \mathcal{A}_n \rightarrow \baxt_3\times \baxt_3$ by
\[
k \mapsto (\eta_1(k), \eta_2(k)).
\]
This map can be extended to a homomorphism from $\mathcal{A}_n^{\star}$ to $\baxt_3\times \baxt_3$. Further, $\eta_{ij}$ is also a homomorphism from $(\mathcal{A}_n^{\star},~^{\sharp})$ to $(\baxt_3\times \baxt_3,~^{\sharp})$. This is because for any $k \in \mathcal{A}_n$, $\eta_1(k^{\sharp})=(\eta_2(k))^{\sharp}, (\eta_1(k))^{\sharp}=\eta_2(k^{\sharp})$ which follows from
\begin{align*}
\begin{cases}
\eta_1(k^{\sharp})=(\eta_2(k))^{\sharp}=[\varepsilon]_{\baxt_3},\;\; \;(\eta_1(k))^{\sharp}= \eta_2(k^{\sharp})=[3]_{\baxt_3} & \text{if}\ k = i_1,\\
\eta_1(k^{\sharp})=(\eta_2(k))^{\sharp}=[\varepsilon]_{\baxt_3},\;\;\; (\eta_1(k))^{\sharp}= \eta_2(k^{\sharp})=[32]_{\baxt_3} & \text{if}\ i_1 < k < i_2,\\
\eta_1(k^{\sharp})=(\eta_2(k))^{\sharp}=[\varepsilon]_{\baxt_3},\;\;\; (\eta_1(k))^{\sharp}= \eta_2(k^{\sharp})=[2]_{\baxt_3}  & \text{if}\ k= i_2,\\
\eta_1(k^{\sharp})=(\eta_2(k))^{\sharp}=[2]_{\baxt_3},\;\;\; (\eta_1(k))^{\sharp}= \eta_2(k^{\sharp})=[\varepsilon]_{\baxt_3}  & \text{if}\ k= i_3,\\
\eta_1(k^{\sharp})=(\eta_2(k))^{\sharp}=[21]_{\baxt_3},\; (\eta_1(k))^{\sharp}= \eta_2(k^{\sharp})=[\varepsilon]_{\baxt_3}& \text{if}\  i_3< k< i_4,\\
\eta_1(k^{\sharp})=(\eta_2(k))^{\sharp}=[1]_{\baxt_3},\;\;\; (\eta_1(k))^{\sharp}=\eta_2(k^{\sharp})=[\varepsilon]_{\baxt_3} & \text{if}\ k = i_4,\\
\eta_1(k^{\sharp})=(\eta_2(k))^{\sharp}=[\varepsilon]_{\baxt_3},\;\;\; (\eta_1(k))^{\sharp}= \eta_2(k^{\sharp})=[\varepsilon]_{\baxt_3} & \text{otherwise}.
\end{cases}
\end{align*}
Therefore it is routine to verify that $\eta_{ij}(\bfw^{\sharp})=(\eta_{ij}(\bfw))^{\sharp}$ for any $\bfw\in \mathcal{A}_n^{\star}$.
%\begin{align*}
%\eta_{ij}(w^{\sharp})&= (\eta_1(w^{\sharp}), \eta_2(w^{\sharp}))\\
%&=(\eta_1(k_n^{\sharp})\cdots\eta_1(k_1^{\sharp}), \eta_2(k_n^{\sharp})\cdots\eta_2(k_1^{\sharp}))\\
%&=( (\eta_2(k_n))^{\sharp}\cdots(\eta_2(k_1))^{\sharp}, (\eta_1(k_n))^{\sharp}\cdots(\eta_1(k_1))^{\sharp})\\
%&=((\eta_2(w))^{\sharp}, (\eta_1(w))^{\sharp})\\
%&=(\eta_{ij}(w))^{\sharp}.
%\end{align*}
%Note that $\eta_{ij}(w)$ is a pair of baxtplactic class of the word obtained from $w$. The first component can be obtained from $w$ by replacing any occurrence of $i$ by $1$; any occurrence of $j$ by $2$; any occurrence of $k$, with $i < k < j$, by $21$; and erasing any occurrence of any other element.
%The second component can be obtained from $w$ by replacing any occurrence of an $j^{\sharp}$ by $2$; any occurrence of $i^{\sharp}$ by $3$; any occurrence of $k$, with $j^{\sharp} < k < i^{\sharp}$, by $32$; and erasing any occurrence of any other element.

{\bf Case 4.} $i<j^{\sharp}<j<i^{\sharp}$ or $j^{\sharp}<i<i^{\sharp}<j$. For convenience, let $i_1=i, i_2=j^{\sharp}, i_3=j$ and $i_4=i^{\sharp}$ when $i<j^{\sharp}<j<i^{\sharp}$ and $i_1=j^{\sharp}, i_2=i,i_3=i^{\sharp}$ and $i_4=j$ when $j^{\sharp}<i<i^{\sharp}<j$. Define  maps $\kappa_1: \mathcal{A}_n \rightarrow \baxt_3$ and $\kappa_2: \mathcal{A}_n \rightarrow \baxt_3$ by
\begin{align*}
k \mapsto \begin{cases}
[1]_{\baxt_3} & \text{if}\ k = i_1,\\
[21]_{\baxt_3} & \text{if}\ i_1 < k < i_3,\\
[2]_{\baxt_3} & \text{if}\  k= i_3,\\
\left[\varepsilon\right]_{\baxt_3} & \text{otherwise},
\end{cases} \;\;\;
k \mapsto \begin{cases}
[2]_{\baxt_3} & \text{if}\  k = i_2,\\
[3]_{\baxt_3} & \text{if}\ k = i_4,\\
[32]_{\baxt_3} & \text{if}\ i_2 < k < i_4,\\
\left[\varepsilon\right]_{\baxt_3} & \text{otherwise}
\end{cases}
\end{align*}
respectively.
Clearly, $\kappa_1, \kappa_2$ can be extended to homomorphisms from $\mathcal{A}_n^{\star}$ to $\baxt_3$ respectively.
Define a map $\kappa_{ij}: \mathcal{A}_n \rightarrow \baxt_3\times \baxt_3$,
\[
k \mapsto (\kappa_1(k), \kappa_2(k)).
\]
Clearly, this map can be extended to a homomorphism from $\mathcal{A}_n^{\star}$ to $\baxt_3\times \baxt_3$. Further, $\kappa_{ij}$ is a homomorphism from $(\mathcal{A}_n^{\star},~^{\sharp})$ to $(\baxt_3\times \baxt_3,~^{\sharp})$.
This is because for any $k \in \mathcal{A}_n$, $\kappa_1(k^{\sharp})=(\kappa_2(k))^{\sharp}, (\kappa_1(k))^{\sharp}=\kappa_2(k^{\sharp})$ which follows from
\begin{align*}
\begin{cases}
\kappa_1(k^{\sharp})=(\kappa_2(k))^{\sharp}=[\varepsilon]_{\baxt_3}, \;\;\;(\kappa_1(k))^{\sharp}= \kappa_2(k^{\sharp})=[3]_{\baxt_3} & \text{if}\ k = i_1,\\
\kappa_1(k^{\sharp})=(\kappa_2(k))^{\sharp}=[\varepsilon]_{\baxt_3}, \;\;\;(\kappa_1(k))^{\sharp}=\kappa_2(k^{\sharp})=[32]_{\baxt_3} & \text{if}\ i_1 < k < i_2,\\
\kappa_1(k^{\sharp})=(\kappa_2(k))^{\sharp}=[2]_{\baxt_3},\;\;\; (\kappa_1(k))^{\sharp}=\kappa_2(k^{\sharp})=[32]_{\baxt_3}  & \text{if}\ k = i_2,\\
\kappa_1(k^{\sharp})=(\kappa_2(k))^{\sharp}=[21]_{\baxt_3}, \;(\kappa_1(k))^{\sharp}= \kappa_2(k^{\sharp})=[32]_{\baxt_3}  & \text{if}\ i_2< k< i_3,\\
\kappa_1(k^{\sharp})=(\kappa_2(k))^{\sharp}=[21]_{\baxt_3},\;  (\kappa_1(k))^{\sharp}= \kappa_2(k^{\sharp})=[2]_{\baxt_3}  & \text{if}\ k = i_3,\\
\kappa_1(k^{\sharp})=(\kappa_2(k))^{\sharp}=[21]_{\baxt_3}, \;(\kappa_1(k))^{\sharp}= \kappa_2(k^{\sharp})=[\varepsilon]_{\baxt_3}& \text{if}\  i_3< k< i_4,\\
\kappa_1(k^{\sharp})=(\kappa_2(k))^{\sharp}=[1]_{\baxt_3},\;\;\; (\kappa_1(k))^{\sharp}=\kappa_2(k^{\sharp})=[\varepsilon]_{\baxt_3} & \text{if}\ k=i_4,\\
\kappa_1(k^{\sharp})=(\kappa_2(k))^{\sharp}=[\varepsilon]_{\baxt_3},\;\;\;(\kappa_1(k))^{\sharp}= \kappa_2(k^{\sharp})=[\varepsilon]_{\baxt_3} & \text{otherwise}.
\end{cases}
\end{align*}
Therefore it is routine to verify that $\kappa_{ij}(\bfw^{\sharp})=(\kappa_{ij}(\bfw))^{\sharp}$ for any $\bfw\in \mathcal{A}_n^{\star}$.
%\begin{align*}
%\kappa_{ij}(w^{\sharp})&= (\kappa_1(w^{\sharp}), \kappa_2(w^{\sharp}))\\
%&=(\kappa_1(w_n^{\sharp})\cdots\kappa_1(w_1^{\sharp}), \kappa_2(w_n^{\sharp})\cdots\kappa_2(w_1^{\sharp}))\\
%&=( (\kappa_2(w_n))^{\sharp}\cdots(\kappa_2(w_1))^{\sharp}, (\kappa_1(w_n))^{\sharp}\cdots(\kappa_2(w_1))^{\sharp})\\
%&=((\kappa_2(w))^{\sharp}, (\kappa_1(w))^{\sharp})\\
%&=(\kappa_{ij}(w))^{\sharp}.
%\end{align*}
%
%Note that $\kappa_{ij}(w)$ is a pair of baxtplactic class of the word obtained from $w$. The first component can be obtained from $w$ by replacing any occurrence of $i$ by $1$; any occurrence of $k$, with $i < k < j$, by $21$; any occurrence of an $j$ by $2$; and erasing any occurrence of any other element. The second component can be obtained from $w$ by replacing any occurrence of an $j^{\sharp}$ by $2$; any occurrence of $i^{\sharp}$ by $3$; any occurrence of $k$, with $j^{\sharp} < k < i^{\sharp}$, by $32$; and erasing any occurrence of any other element.

Now we can define the map $\varphi_{ij}: \mathcal{A}^{\star}_n \rightarrow \baxt_3\times \baxt_3$ by
\begin{align*}
\varphi_{ij}=
\left\{
\begin{array}{ll}
\lambda_{ij} & \hbox{if $i^{\sharp}=j$,} \\
[0.1cm]
\theta_{ij} & \hbox{if $i<j=j^{\sharp}<i^{\sharp}$ or $j^{\sharp}<i=i^{\sharp}<j$,} \\
[0.1cm]
\eta_{ij} & \hbox{if $i<j<j^{\sharp}<i^{\sharp}$ or $j^{\sharp}<i^{\sharp}<i<j$,}\\
[0.1cm]
\kappa_{ij} & \hbox{if $i<j^{\sharp}<j<i^{\sharp}$ or $j^{\sharp}<i<i^{\sharp}<j$}\\
\end{array}
\right.
\end{align*}
where $\lambda_{ij}, \theta_{ij}, \eta_{ij}, \kappa_{ij}$ are defined as above. Since each of maps $\lambda_{ij}, \theta_{ij}, \eta_{ij}, \kappa_{ij}$ is a homomorphism from $(\mathcal{A}_n^{\star},~^{\sharp})$ to $(\baxt_3\times \baxt_3,~^{\sharp})$.
Therefore $\varphi_{ij}$ is a homomorphism from $(\mathcal{A}_n^{\star},~^{\sharp})$ to $(\baxt_3\times \baxt_3,~^{\sharp})$.

\begin{lemma}\label{lemma:phi_factor_homomorphism}
The homomorphism $\varphi_{ij}$ induces a homomorphism $\varphi_{ij} : (\baxt_n,~^{\sharp}) \rightarrow (\baxt_3\times \baxt_3,~^{\sharp})$ for any $n\geq 4$.
\end{lemma}
	
\begin{proof}
Note that $\varphi_{ij}$ is a homomorphism from $(\mathcal{A}_n^{\star}, ~^{\sharp})$ to $(\baxt_3\times \baxt_3,~^{\sharp})$. Then to show that the homomorphism $\varphi_{ij}$ induces a homomorphism from $(\mathsf{baxt}_n,~^\sharp)$ to $(\baxt_3\times \baxt_3,~^{\sharp})$, we only need to show that for any $\bfu,\bfv\in \mathcal{A}_n^{\star}$, if $\bfu\equiv_{\baxt_\infty}\bfv$, then $\varphi_{ij}(\bfu)=\varphi_{ij}(\bfv)$.
It follows from Proposition \ref{pro:inversion} and the definition of $\varphi_{ij}$ that the evaluations of the first and the second components of $\varphi_{ij}(\bfu)$ and $\varphi_{ij}(\bfv)$ are the same respectively. In the following, we prove that the left and right precedences of the first and the second components of $\varphi_{ij}(\bfu)$ and $\varphi_{ij}(\bfv)$ are the same respectively.

{\bf Case 1.} $\varphi_{ij}=\lambda_{ij}$. If there exists $k\in \supp(\bfu)$ satisfying $i<k<j$, since the $1$-$3$  left precedence in the first component of $\lambda_{ij}(\bfu)$ [resp. $\lambda_{ij}(\bfv)$] corresponds to some $i$-$h$ left  precedence in $\bfu$ [resp. $\bfv$] with $i<h\leq j$ and the $3$-$1$  right precedence in the first component of $\lambda_{ij}(\bfu)$ [resp. $\lambda_{ij}(\bfv)$] corresponds to some $h$-$i$ right  precedence in $\bfu$ [resp. $\bfv$] with $i<h\leq j$, it follows from Proposition \ref{pro:inversion} that the first components of $\lambda_{ij}(\bfu)$ and $\lambda_{ij}(\bfv)$ have the same  $1$-$3$ left  precedence and the same $3$-$1$ right precedence; if there is no $k\in\supp(\bfu)$ satisfying $i<k<j$, then since the $1$-$3$  left precedence in the first component of $\lambda_{ij}(\bfu)$ [resp. $\lambda_{ij}(\bfv)$] corresponds to the $i$-$j$ left  precedence in $\bfu$ [resp. $\bfv$] and the $3$-$1$  right precedence in the first component of $\lambda_{ij}(\bfu)$ [resp. $\lambda_{ij}(\bfv)$] corresponds to the $j$-$i$ right  precedence in $\bfu$ [resp. $\bfv$], it follows from  Proposition \ref{pro:inversion} that the first component of $\lambda_{ij}(\bfu)$ and the first component of $\lambda_{ij}(\bfv)$ have the same  $1$-$3$ left  precedence and the same $3$-$1$ right precedence.  A similar argument can show that the left and right precedences of the second components of $\varphi_{ij}(\bfu)$ and $\varphi_{ij}(\bfv)$ are the same.

{\bf Case 2.} $\varphi_{ij}=\theta_{ij}, \eta_{ij}$ or $\kappa_{ij}$.

{\bf 2.1.} $i<j=j^{\sharp}<i^{\sharp}$, $i<j<j^{\sharp}<i^{\sharp}$ or $i<j^{\sharp}<j<i^{\sharp}$.
If there exists $k\in \supp(\bfu)$ satisfying $i<k<j$, since the $1$-$2$  left precedence in the first component of $\varphi_{ij}(\bfu)$ [resp. $\varphi_{ij}(\bfv)$] corresponds to some $i$-$h$ left  precedence in $\bfu$ [resp. $\bfv$] with $i<h\leq j$ and the $2$-$1$  right precedence in the first component of $\varphi_{ij}(\bfu)$ [resp. $\varphi_{ij}(\bfv)$] corresponds to some $h$-$i$ right  precedence in $\bfu$ [resp. $\bfv$] for $i<h\leq j$, it follows from Proposition \ref{pro:inversion} that the first components of $\varphi_{ij}(\bfu)$ and $\varphi_{ij}(\bfv)$ have the same  $1$-$2$ left  precedence and the same $2$-$1$ right precedence; if there is no $k\in \supp(\bfu)$ satisfying $i<k<j$, then since the $1$-$2$ left  precedence in the first component of $\varphi_{ij}(\bfu)$ [resp. $\varphi_{ij}(\bfv)$] corresponds to the $i$-$j$ left  precedence in $\bfu$ [resp. $\bfv$] and the $2$-$1$ right  precedence in the first component of $\varphi_{ij}(\bfu)$ [resp. $\varphi_{ij}(\bfv)$] corresponds to the $j$-$i$ right  precedence in $\bfu$ [resp. $\bfv$], it follows from Proposition \ref{pro:inversion} that the first component of $\varphi_{ij}(\bfu)$ and the first component of $\varphi_{ij}(\bfv)$ have the same  $1$-$2$ left  precedence and the same $2$-$1$ right precedence. A similar argument can show that the left and right precedences of the second components of $\varphi_{ij}(\bfu)$ and $\varphi_{ij}(\bfv)$ are the same.

{\bf 2.2.}  $j^{\sharp}<i=i^{\sharp}<j$, $j^{\sharp}<i^{\sharp}<i<j$, or $j^{\sharp}<i<i^{\sharp}<j$. If there exists $k\in \supp(\bfu)$ satisfying $i<k<j$, since the $2$-$3$  left precedence in the second component of $\varphi_{ij}(\bfu)$ [resp. $\varphi_{ij}(\bfv)$] corresponds to some $i$-$h$ left  precedence in $\bfu$ [resp. $\bfv$] with $i<h\leq j$ and the $3$-$2$  right precedence in the second component of $\varphi_{ij}(\bfu)$ [resp. $\varphi_{ij}(\bfv)$] corresponds to some $h$-$i$ right  precedence in $\bfu$ [resp. $\bfv$] with $i<h\leq j$, it follows from Proposition \ref{pro:inversion} that the second components of $\varphi_{ij}(\bfu)$ and $\varphi_{ij}(\bfv)$ have the same  $2$-$3$ left  precedence and the same $3$-$2$ right precedence; if there is no $k\in \supp(\bfu)$ satisfying $i<k<j$, then since the $2$-$3$ left  precedence in the second component of $\varphi_{ij}(\bfu)$ [resp. $\varphi_{ij}(\bfv)$] corresponds to some $i$-$j$ left  precedence in $\bfu$ [resp. $\bfv$] and the $3$-$2$ right  precedence in the second component of $\varphi_{ij}(\bfu)$ [resp. $\varphi_{ij}(\bfv)$] corresponds to the $j$-$i$ right  precedence in $\bfu$ [resp. $\bfv$], it follows from Proposition \ref{pro:inversion} that the second component of $\varphi_{ij}(\bfu)$ and the second component of $\varphi_{ij}(\bfv)$ have the same  $2$-$3$ left  precedence and the same $3$-$2$ right precedence. A similar argument can show that the left and right precedences of the first components of $\varphi_{ij}(\bfu)$ and $\varphi_{ij}(\bfv)$ are the same.
\end{proof}
	
\begin{lemma}\label{coro:varphi_equiv_baxtn}
Let $\bfu, \bfv \in \mathcal{A}_n^{\star}$ for any $n\geq 4$. Then $\bfu \equiv_{\baxt_\infty} \bfv$ if and only if $\varphi_{ij}(\bfu) = \varphi_{ij}(\bfv)$ for all $1 \leq i < j \leq n$.
\end{lemma}

\begin{proof}
The necessity follows from the proof of Lemma~\ref{lemma:phi_factor_homomorphism}. Let $\bfw \in \mathcal{A}_n^{\star}$ for some $n \geq 4$. Suppose $\supp(\bfw)= \{a_1 < \dots < a_{\ell}\}$ for some $\ell \in \mathbb{N}$. We can obtain the evaluation of $\bfw$ from $\varphi_{a_i, a_j}(\bfw)$. For any $i$ with $1 \leq i < \ell$, if $a_i^{\sharp}=a_{i+1}$, then since there is no $k\in\supp(\bfw)$ satisfying $a_i<k<a_{i+1}$, it follows from the definition of $\varphi_{ij}$ that the number of occurrences of $1$ in the first [resp. second] component of $\varphi_{ij}(\bfw)$ equals to ${|\bfw|}_{a_i}$ and the number of occurrences of $3$ in the first [resp. second] component of $\varphi_{ij}(\bfw)$ equals to ${|\bfw|}_{a_{i+1}}$; if
$a_i<a_{i+1}=a_{i+1}^{\sharp}<a_i^{\sharp}$, $a_i<a_{i+1}<a_{i+1}^{\sharp}<a_i^{\sharp}$ or $a_i<a_{i+1}^{\sharp}<a_{i+1}<a_i^{\sharp}$, then since there is no $k\in\supp(\bfw)$ satisfying $a_i<k<a_{i+1}$, it follows from the definition of $\varphi_{ij}$ that the the number of occurrences of $1$ in the first component of $\varphi_{ij}(\bfw)$ equals to ${|\bfw|}_{a_i}$ and the number of occurrences of $2$ in the first component of $\varphi_{ij}(\bfw)$ equals to ${|\bfw|}_{a_{i+1}}$; if $a_{i+1}^{\sharp}<a_{i}=a_{i}^{\sharp}<a_{i+1}$, $a_{i+1}^{\sharp}<a_{i}^{\sharp}<a_{i}<a_{i+1}$, or $a_{i+1}^{\sharp}<a_{i}<a_{i}^{\sharp}<a_{i+1}$, then since there is no $k\in\supp(\bfw)$ satisfying $a_i<k<a_{i+1}$, it follows from the definition of $\varphi_{ij}$ that the number of occurrences of $2$ in the second component of $\varphi_{ij}(\bfw)$ equals to ${|\bfw|}_{a_i}$ and the number of occurrences of $3$ in the second component of $\varphi_{ij}(\bfw)$ equals to ${|\bfw|}_{a_{i+1}}$.
Therefore both the number of occurrences of $a_i$ and $a_{i+1}$ in $\bfw$ can be derived from the maps $\varphi_{a_i a_{i+1}}$.

We can obtain the left and right precedences of $\bfw$ from $\varphi_{a_i, a_j}(\bfw)$. If $a_i,a_j$ satisfy $a_i^{\sharp}=a_j$, $a_i<a_{j}=a_{j}^{\sharp}<a_i^{\sharp}$, $a_i<a_{j}<a_{j}^{\sharp}<a_i^{\sharp}$ or $a_i<a_{j}^{\sharp}<a_{j}<a_i^{\sharp}$, then we check the the first component of $\varphi_{a_i, a_j}(\bfw)$; if $a_i,a_j$ satisfy $a_{j}^{\sharp}<a_{i}=a_{i}^{\sharp}<a_{j}$, $a_{j}^{\sharp}<a_{i}^{\sharp}<a_{i}<a_{j}$, or $a_{j}^{\sharp}<a_{i}<a_{i}^{\sharp}<a_{j}$, then we check the the second component of $\varphi_{a_i, a_j}(\bfw)$.

For each $a_j$ satisfying $a_{j}^{\sharp}<a_{j}$, if $a_i$ satisfying $a_{j}^{\sharp}<a_{i}<a_{j}$ is the largest number such that the second component of $\varphi_{a_i,a_j}(\bfw)$ starts with $2$,
then, when reading the second component of $\varphi_{a_i,a_j}(\bfw)$ from left-to-right, the first occurrence of $3$  corresponds to the first occurrence of $a_j$,  and all
occurrences of $2$ before the first occurrence of $3$ correspond to all occurrences of $a_i$ before the first occurrence of $a_j$. Hence $\bfw$ has a $a_i$-$a_j$ left precedence. Otherwise, if $a_i$ satisfying $a_i^{\sharp}=a_j$ or $a_i<a_{j}^{\sharp}<a_{j}<a_i^{\sharp}$ is the largest number such that the first component of $\varphi_{a_i,a_j}(\bfw)$ starts with $1$,
then, when reading the first component of $\varphi_{a_i,a_j}(\bfw)$ from left-to-right, the first occurrence of $2$ or $3$  corresponds to the first occurrence of $a_j$, and all
occurrences of $1$ before the first occurrence of $2$ or $3$ correspond to all occurrences of $a_i$ before the first occurrence of $a_j$, and so $\bfw$ has a $a_i$-$a_j$ left precedence; otherwise $\bfw$ does not have a $a_i$-$a_j$ left precedence.
For each $a_j$ satisfying $ a_j\leq a_j^{\sharp}$, if $a_i<a_j$ is the largest number such that the first component of $\varphi_{a_i,a_j}(\bfw)$ starts with $1$,
then, when reading the first component of $\varphi_{a_i,a_j}(\bfw)$ from left-to-right, the first occurrence of $2$  corresponds to the first occurrence of $a_j$, and all
occurrences of $1$ before the first occurrence of $2$ correspond to all occurrences of $a_i$ before the first occurrence of $a_j$, and so $\bfw$ has a $a_i$-$a_j$ left precedence; otherwise $\bfw$ does not have a $a_i$-$a_j$ left precedence.

For each $a_i$  satisfying $a_i<a_i^{\sharp}$, if $a_j$ satisfying $a_{i}<a_{j}\leq a_i^{\sharp}$ is the smallest number such that the first component of $\varphi_{a_i,a_j}(\bfw)$ ends with $2$ or $3$,
then, when reading the first component of $\varphi_{a_i,a_j}(\bfw)$ from right-to-left, the first occurrence of $1$  corresponds to the first occurrence of $a_i$,  and all
occurrences of $2$ or $3$ before the first occurrence of $1$ correspond to all occurrences of $a_j$ before the first occurrence of $a_i$. Hence $\bfw$ has a $a_j$-$a_i$ right precedence. Otherwise, if  $a_j$ satisfying $a_{j}^{\sharp}<a_{i}<a_{i}^{\sharp}<a_{j}$ is the smallest number such that the second component of $\varphi_{a_i,a_j}(\bfw)$ ends with $3$,
then, when reading the second component of $\varphi_{a_i,a_j}(\bfw)$ from right-to-left, the first occurrence of $2$  corresponds to the first occurrence of $a_i$, and all
occurrences of $3$ before the first occurrence of $2$ correspond to all occurrences of $a_j$ before the first occurrence of $a_i$, and so $\bfw$ has a $a_j$-$a_i$ right precedence; otherwise $\bfw$ does not have a $a_j$-$a_i$ right precedence.
For each $a_i$ satisfying $a_i\geq a_i^{\sharp}$, if $a_j$ is the smallest number such that the second component of $\varphi_{a_i,a_j}(\bfw)$ ends with $3$,
then, when reading the second component of $\varphi_{a_i,a_j}(\bfw)$ from right-to-left, the first occurrence of $2$ corresponds to the first occurrence of $a_i$, and all
occurrences of $3$ before the first occurrence of $2$ correspond to all occurrences of $a_j$ before the first occurrence of $a_i$, and so $\bfw$ has a $a_j$-$a_i$ right precedence; otherwise $\bfw$ does not have a $a_j$-$a_i$ right precedence.

Therefore, if $\varphi_{ij}(\bfu) = \varphi_{ij}(\bfv)$ for all $1 \leq i < j \leq n$, then by the above arguments, $\bfu$ and $\bfv$ have the same evaluations and left and right precedences. Consequently, $\bfu \equiv_{\baxt_{\infty}} \bfv$.
\end{proof}

For each $n \in \mathbb{N}$, with $n \geq 4$, let $I_n$ be the index set
\[
\{ (i,j): 1 \leq i < j \leq n \}.
\]
Now, consider the map
\[
\phi_n : (\baxt_n,~^{\sharp}) \rightarrow \prod\limits_{I_n} (\baxt_3\times \baxt_3,~^{\sharp}),
\]
whose $(i,j)$-th component is given by $\varphi_{ij}([\bfw]_{\baxt_n})$ for $\bfw \in \mathcal{A}_n^{\star}$ and $(i,j) \in I_n$.
		
\begin{proposition}\label{prop:baxtn_embed_baxt3}
The map $\phi_n$ is an embedding from $(\baxt_n,~^{\sharp})$ to $(\baxt_3\times \baxt_3,~^{\sharp})$.
\end{proposition}
	
\begin{proof}
The map $\phi_n$ is a homomorphism by Lemma~\ref{lemma:phi_factor_homomorphism}. It follows from the definition of $\phi_n$ and Lemma~\ref{coro:varphi_equiv_baxtn} that $[\bfu]_{\baxt_n} =[\bfv]_{\baxt_n}$ if and only if $\varphi_n([\bfu]_{\baxt_n}) = \varphi_n([\bfv]_{\baxt_n})$ for any $\bfu, \bfv \in \mathcal{A}_n^{\star}$. Hence $\phi_n$ is an embedding.
\end{proof}

For any $([\bfu_1]_{\baxt_3}, [\bfu_2]_{\baxt_3}, \dots, [\bfu_{2|I_n|}]_{\baxt_3})\in \prod\limits_{2I_n}\baxt_3$, define an involution operation $^{\sharp}$ on $\prod\limits_{2I_n}\baxt_3$ by
\begin{align*}
([\bfu_1]_{\baxt_3}, [\bfu_2]_{\baxt_3}, \dots, [\bfu_{2|I_n|}]_{\baxt_3})^{\sharp}=([\bfu_{2|I_n|}]_{\baxt_3}^{\sharp}, \dots, [\bfu_{2}]_{\baxt_3}^{\sharp}, [\bfu_1]_{\baxt_3}^{\sharp}).
\end{align*}
Define a map $\psi:\prod\limits_{I_n} (\baxt_3\times \baxt_3,~^{\sharp})\rightarrow (\prod\limits_{2I_n}\baxt_3,~^{\sharp})$ given by
\begin{align*}
(([\bfu_1]_{\baxt_3}, [\bfu_2]_{\baxt_3}), ([\bfu_3]_{\baxt_3}, [\bfu_2]_{\baxt_4}), \dots, ([\bfu_{2|I_n|-1}]_{\baxt_3}, [\bfu_{2|I_n|}]_{\baxt_3}))\\
\mapsto ([\bfu_1]_{\baxt_3}, [\bfu_3]_{\baxt_3}, \dots, [\bfu_{2|I_n|-1}]_{\baxt_3}, [\bfu_{2|I_n|}]_{\baxt_3}, \dots, [\bfu_4]_{\baxt_3}, [\bfu_2]_{\baxt_3}).
\end{align*}
It is routine to verify that the map $\psi$ is an isomorphism. By Theorem \ref{thm:baxt3-repre}, each element in $(\baxt_3,~^{\sharp})$ corresponds to a matrix in $UT_{15}(\mathbb{S})$ and the involution $^{\sharp}$ on $(\baxt_3,~^{\sharp})$ corresponds to the skew transposition on $UT_{15}(\mathbb{S})$. It follows that there is an embedding, denoted by $\varphi$, from $(\prod\limits_{2I_n}\baxt_3,~^{\sharp})$ to $(UT_{30|I_n|}(\mathbb{S}),~^{D})$. Let
\[
\varphi_n=\varphi\circ\psi\circ\phi_n.
\]
Then the following result hold.

\begin{theorem}\label{thm_baxt_tropical}
For each $n\geq 4$,  the map $\varphi_n:(\mathsf{baxt}_n,~^\sharp)\rightarrow(UT_{30|I_n|}(\mathbb{S}),~^{D})$ is a faithful representation of $(\mathsf{baxt}_n,~^\sharp)$.
\end{theorem}

Let $a<b<c<d$ be a 4-element ordered alphabet and
\[
B=\langle a,b,c,d\,|\,\mathcal{R}_{\baxt_\infty}, ac=ca, ad=da, bc=cb, bd=db\rangle\cup \{1\}
\]
be a monoid.
The involution operation $^{\sharp}$ on $B$ can be defined by $a\mapsto d, b\mapsto c$.

\begin{theorem}\label{thm:same-var}
For any $m,n\geq 4$, the involution monoids $(\baxt_m,~^{\sharp})$ and $(\baxt_n,~^{\sharp})$ generate the same variety.
\end{theorem}

\begin{proof}
It follows from Proposition \ref{prop:baxtn_embed_baxt3} that
\[
\var(\baxt_4,~^{\sharp})\subseteq \var(\baxt_5,~^{\sharp}) \subseteq \cdots \subseteq\var (\baxt_3\times \baxt_3,~^{\sharp}).
\]
It suffices to show that $\var (\baxt_3\times \baxt_3,~^{\sharp})\subseteq \var(\baxt_4,~^{\sharp})$.
Clearly $(B, ~^{\sharp}$) is a homomorphism image of $(\baxt_4,~^{\sharp})$, and so $\var(B,~^{\sharp})\subseteq \var(\baxt_4,~^{\sharp})$. In the following, we show that $\var (\baxt_3\times \baxt_3,~^{\sharp})\subseteq \var(B,~^{\sharp})$.

Let $\varphi$ be a map from $(\baxt_3\times \baxt_3,~^{\sharp})$ to $(B,~^{\sharp})\times (B,~^{\sharp})\times (B,~^{\sharp})$ given by
\begin{align*}
&([\varepsilon]_{\baxt_3}, [\varepsilon]_{\baxt_3})\mapsto (1,1,1),\\
([1]_{\baxt_3},[\varepsilon]_{\baxt_3})&\mapsto (a,a,1), \quad\quad([2]_{\baxt_3},[\varepsilon]_{\baxt_3})\mapsto (b, ba,a), \\ ([3]_{\baxt_3},[\varepsilon]_{\baxt_3})&\mapsto (1,b,b), \quad\quad\;([\varepsilon]_{\baxt_3}, [1]_{\baxt_3})\mapsto (1,c,c), \\
([\varepsilon]_{\baxt_3}, [2]_{\baxt_3})&\mapsto (c,dc,d), \quad\;\;\;([\varepsilon]_{\baxt_3}, [3]_{\baxt_3})\mapsto (d,d,1).
\end{align*}
First to show that $\varphi$ is well-defined, that is, if $\bfu_1\equiv_{\baxt_{\infty}}\bfv_1, \bfu_2\equiv_{\baxt_{\infty}}\bfv_2$, then
\[
(U_1,U_2,U_3)=\varphi([\bfu_1]_{\baxt_3},[\bfu_2]_{\baxt_3})=\varphi([\bfv_1]_{\baxt_3},[\bfv_2]_{\baxt_3})=(V_1,V_2,V_3).
\]
Since $\bfu_1\equiv_{\baxt_{\infty}}\bfv_1, \bfu_2\equiv_{\baxt_{\infty}}\bfv_2$, it is easy to see that $\ev(U_i)=\ev(V_i)$ for $i=1,2,3$.
It follows from Proposition \ref{pro:inversion} and the definition of $(B, ~^{\sharp})$ that, for any $([\bfw_1]_{\baxt_3}, [\bfw_2]_{\baxt_3})\in\baxt_3\times \baxt_3$, the left and right precedences of each component of $\varphi([\bfw_1]_{\baxt_3}, [\bfw_2]_{\baxt_3})=(W_1, W_2, W_3)$  can be characterized as follows:
\begin{align*}
(a\textrm{-}b, \ell_1) \in\lpi(W_1)~[\text{resp. } (c\textrm{-}d, \ell_1) \in\lpi(W_3)]&\Leftrightarrow  (1\textrm{-}2, \ell_1)\in \lpi(\bfw_1)~[\text{resp. } \lpi(\bfw_2)],\\
(a\textrm{-}b, \ell_2) \in\lpi(W_3)~[\text{resp. } (c\textrm{-}d, \ell_2) \in\lpi(W_1)]&\Leftrightarrow (2\textrm{-}3, \ell_2) \in \lpi(\bfw_1)~[\text{resp. } \lpi(\bfw_2)],\\
(a\textrm{-}b, \ell_3)\in\lpi(W_2)~[\text{resp. } (c\textrm{-}d,\ell_3) \in\lpi(W_2)] &\Leftrightarrow (1\textrm{-}3, \ell_3)\in \lpi(\bfw_1)~[\text{resp. } \lpi(\bfw_2)],
\end{align*}
where $\ell_3\neq \ell_1$ and
\begin{align*}
(b\textrm{-}a, r_1) \in\rpi(W_1)~[\text{resp. } (d\textrm{-}c, r_1) \in\rpi(W_3)]&\Leftrightarrow  (2\textrm{-}1, r_1)\in \rpi(\bfw_1)~[\text{resp. } \rpi(\bfw_2)],\\
(b\textrm{-}a, r_2) \in\rpi(W_3)~[\text{resp. } (d\textrm{-}c, r_2) \in\rpi(W_1)]&\Leftrightarrow (3\textrm{-}2, r_2) \in \rpi(\bfw_1)~[\text{resp. } \rpi(\bfw_2)],\\
(b\textrm{-}a, r_3)\in\rpi(W_2)~[\text{resp. } (d\textrm{-}c,r_3) \in\rpi(W_2)] &\Leftrightarrow (3\textrm{-}1, r_3)\in \rpi(\bfw_1)~[\text{resp. } \rpi(\bfw_2)].
\end{align*}
where $r_3\neq r_1$.
Since $\bfu_1\equiv_{\baxt_{\infty}}\bfv_1, \bfu_2\equiv_{\baxt_{\infty}}\bfv_2$, it is routine to show that $\lpi(U_i)=\lpi(V_i)$ and $\rpi(U_i)=\rpi(V_i)$ for $i=1,2,3$. Thus $\varphi$ is well-defined, and so $\varphi$ is a homomorphism.

Next to show that the homomorphism $\varphi$ is injective.
Suppose that $([\bfu_1]_{\baxt_3},$ $[\bfu_2]_{\baxt_3})\neq ([\bfv_1]_{\baxt_3}, [\bfv_2]_{\baxt_3})$. Then $\ev(\bfu_1)\neq\ev(\bfu_2)$, or $\ev(\bfv_1)\neq\ev(\bfv_2)$, or $\lpi(\bfu_1)\neq\lpi(\bfu_2)$, or $\rpi(\bfu_1)\neq\rpi(\bfu_2)$, or $\lpi(\bfv_1)\neq\lpi(\bfv_2)$ or $\rpi(\bfv_1)\neq\rpi(\bfv_2)$.  Hence, by the definition of $\varphi$, it is routine to show that $\varphi([\bfu_1]_{\baxt_3},$ $[\bfu_2]_{\baxt_3})\neq\varphi([\bfv_1]_{\baxt_3},$ $[\bfv_2]_{\baxt_3})$. Therefore the homomorphism $\varphi$ is injective.
\end{proof}

\section{The identities satisfied by $(\mathsf{baxt}_n,~^\sharp)$} \label{sec:characterization} %

%Cain et al. gave a complete characterization of the identities satisfied by $\hypo_n$ for each $n\geq 2$ in \cite[Theorem~4.1]{CMR21a}.
In this section, a complete characterization of the word identities satisfied by $(\mathsf{baxt}_n,~^\sharp)$ for each finite $n$ is given.
Clearly, a word identity $\bfu\approx \bfv$ holds in $(\mathsf{baxt}_1,~^\sharp)$ if and only if $\occ(x,\overline{\bfu})=\occ(x,\overline{\bfv})$ for any $x\in \con(\overline{\bfu\bfv})$.

%Let $\bfu$ be a word and $\mathcal{Y}$ and $\mathcal{Z}$ be two disjoint subsets of $\con(\bfu)$. Then we say that
%\textit{$\mathcal{Y}$ completely  precedes $\mathcal{Z}$} in $\bfu$, denoted by $\mathcal{Y} \prec_{\bfu}\mathcal{Z}$, if the last occurrence of each variable in $\mathcal{Y}$ precedes the first occurrence of any variable in $\mathcal{Z}$ in $\bfu$.
%%In other words, if $\mathcal{Y}$ completely precedes $\mathcal{Z}$, then $\bfu[\mathcal{Y}\cup\mathcal{Z}] \in \mathcal{Y}^{+}\cdot\mathcal{Z}^{+}$.
%In particular, if $\mathcal{Y}$ or $\mathcal{Z}$ is a singleton set, say $\mathcal{Y}=\{y\}$, then we write $y \prec_{\bfu} \mathcal{Z}$ rather than $\{y\} \prec_{\bfu} \mathcal{Z}$.

Let
\[
A=\langle\, a, b\, |\, ab=ba\,\rangle=\{a^mb^n\,|\,m,n\geq 0\}
\]
be a monoid.
The monoid $A$ is an involution monoid $(A, ~^*)$ under the unary operation $^*: a^mb^n\mapsto a^nb^m$. A word identity $\bfu \approx \bfv$ is \textit{balanced} if $\mathsf{occ}(x, \bfu)=\mathsf{occ}(x, \bfv)$ for any $x\in \mathsf{con}(\bfu\bfv)$.

\begin{lemma}\label{lem:balanced}
A word identity $\bfu\approx \bfv$ holds in $(A, ~^*)$ if and only if $\bfu\approx \bfv$ is balanced.
\end{lemma}

\begin{proof}
Suppose that $\bfu\approx \bfv$ is a word identity satisfied by $(A, ~^*)$ such that either $\occ(x,\bfu)\ne\occ(x,\bfv)$ or $\occ(x^*,\bfu)\ne \occ(x^*,\bfv)$ for some $x\in \mathcal{X}$. Let $\varphi$ be a homomorphism from $(\mathcal{X}\cup \mathcal{X}^*)^{+}$ to $(A,~^*)$ that maps $x$ to $a$ and any other variable to $1$.
Then $\varphi(\bfu)=a^{\occ(x,\bfu)}b^{\occ(x^*,\bfu)} \ne a^{\occ(x,\bfv)}b^{\occ(x^*,\bfv)}=\varphi(\bfv)$, a contradiction.

Conversely, if $\bfu\approx \bfv$ is balanced, then $\varphi(\bfu)=\varphi(\bfv)$ for any homomorphism $\varphi$ from $(\mathcal{X}\cup \mathcal{X}^*)^{+}$ to $(A,~^*)$ since $(A,~^*)$ is commutative. Therefore  $\bfu\approx \bfv$ holds in $(A, ~^*)$.
\end{proof}

If $\bfu$ and $\bfv$ are words such that $\occ(x, \bfu) =
\occ(x, \bfv)$  for each variable $x$, then we say that $\bfv$ is a \textit{permutation} of $\bfu$.

\begin{theorem}\label{thm:baxt2}
A word identity $\bfu\approx \bfv$ holds in $(\mathsf{baxt}_2,~^\sharp)$ if and only if $\bfu\approx \bfv$ is balanced and satisfies that for any $x, y \in \con(\bfu)$ with $x, x^*\neq y$,
\begin{enumerate}[\rm(I)]
 % \item \begin{enumerate}[\rm(a)]
%          \item  if $\bfu[x]\in (x^*)^{\alpha}x\{x,x^*\}^{\times}$, then $\bfv[x]\in (x^*)^{\alpha}x\{x,x^*\}^{\times}$,
%
%          \item if $\bfu[x]\in \{x,x^*\}^{\times}x(x^*)^{\alpha}$, then $\bfv[x]\in \{x,x^*\}^{\times}x(x^*)^{\alpha}$;
%        \end{enumerate}

  \item \begin{enumerate}[\rm(a)]
          \item  if $\bfu[x,y]\in x^{\alpha}x^*\{x,x^*,y,y^*\}^{\times}$ for some $\alpha\geq 1$, then $\bfv[x,y]\in x^{\alpha}x^*\{x,x^*,y,y^*\}^{\times}$,

          \item if $\bfu[x,y]\in \{x,x^*,y,y^*\}^{\times}x^*x^{\alpha}$ for some $\alpha\geq 1$, then $\bfv[x,y]\in \{x,x^*,y,y^*\}^{\times}x^*x^{\alpha}$.
        \end{enumerate}

  \item \begin{enumerate}[\rm(a)]
          \item  if $\bfu[x,y]\in y^{\alpha}x\{x,x^*,y,y^*\}^{\times}$ for some $\alpha\geq 1$, then $\bfv[x,y]\in y^{\alpha}x\{x,x^*,y,y^*\}^{\times}$,

          \item if $\bfu[x,y]\in \{x,x^*,y,y^*\}^{\times}xy^{\alpha}$ for some $\alpha\geq 1$, then $\bfv[x,y]\in \{x,x^*,y,y^*\}^{\times}xy^{\alpha}$.
        \end{enumerate}

  \item \begin{enumerate}[\rm(a)]
          \item  if $\bfu[x,y]\in \bfa x\{x,x^*,y,y^*\}^{\times}$ with $\con(\bfa)= \{x^*, y^*\}$, then
          $\bfv[x,y]\in \bfa' x\{x,x^*,y,y^*\}^{\times}$ or
           $\bfv[x,y]\in \bfa' y\{x,x^*,y,y^*\}^{\times}$ where $\bfa'$ is a permutation of $\bfa$,
          \item   if $\bfu[x,y]\in \{x, x^*, y, y^*\}^{\times} x\bfa$ with $\con(\bfa)= \{x^*, y^*\}$, then
          $\bfv[x,y]\in \{x,x^*,y,y^*\}^{\times}x\bfa'$ or
           $\bfv[x,y]\in \{x,x^*,y,y^*\}^{\times} y\bfa'$
            where $\bfa'$ is a permutation of $\bfa$.
        \end{enumerate}
\end{enumerate}
\end{theorem}
\begin{proof}
%First, we prove that a word identity $\bfu\approx \bfv$ satisfied by $(\mathsf{baxt}_2,~^\sharp)$ must be  balanced and satisfy (I)--(III).
Let $\bfu \approx  \bfv$ be any word identity satisfied by $(\mathsf{baxt}_2,~^\sharp)$. It is routine to verify that $(A,~^*)$ is a homomorphic image of $(\mathsf{baxt}_2,~^\sharp)$ under the map given by $[1]_{\mathsf{baxt}_2}\mapsto a$ and $[2]_{\mathsf{baxt}_2}\mapsto b$. Then it follows from Lemma~\ref{lem:balanced} that $\bfu\approx \bfv$ is balanced. Now to show that (Ia)--(IIIa) hold, and then (Ib)--(IIIb) hold by symmetry.  Suppose that $\bfu$ starts with $x$ but $\bfv$ starts with $y \neq x$. If $y=x^*$, then letting $\phi_1$ be the substitution such that $x\mapsto [2]_{{\baxt}_2}, z\mapsto [\varepsilon]_{{\baxt}_2}$ for any $z\neq x, x^*$, we obtain
\begin{align*}
[2]_{\baxt_2}\cdot [\bfs]_{\baxt_2}=\phi_1(\bfu)\neq \phi_1(\bfv)=[1]_{\baxt_2}\cdot [\bft]_{\baxt_2}
\end{align*}
where $\bfs, \bft \in \mathcal{A}_2^{\star}$, a contradiction; if $y\neq x^*$, then letting $\phi_2$ be the substitution such that $x\mapsto [2]_{{\baxt}_2}, y\mapsto [1]_{{\baxt}_2}, z\mapsto [\varepsilon]_{{\baxt}_2}$ for any $z\neq x, y$, we obtain
\begin{align*}
[2]_{\baxt_2}\cdot [\bfs]_{\baxt_2}=\phi_2(\bfu)\neq \phi_2(\bfv)=[1]_{\baxt_2}\cdot [\bft]_{\baxt_2}
\end{align*}
where $\bfs, \bft \in \mathcal{A}_2^{\star}$, a contradiction. Thus $\bfu$ and $\bfv$ start with the same variable.

%It is routine to show that $(B, ~^*)$ is a homomorphic image of $(\mathsf{baxt}_2,~^\sharp)$ under the map given by $[1\cdots 2]_{\mathsf{baxt}_2}\mapsto a$, $[2\cdots 1]_{\mathsf{baxt}_2}\mapsto b$, $[1\cdots 1]_{\mathsf{baxt}_2}\mapsto ab$ and $[2\cdots 2]_{\mathsf{baxt}_2}\mapsto ba$. Therefore Condition (I) holds by Lemma~\ref{lem:ip-fp}.

%If $\bfu[x]\in (x^*)^{\alpha}x\{x,x^*\}^{\times}$,
%then $\bfv[x]\in (x^*)^{\beta}x\{x,x^*\}^{\times}$.
%Suppose that $\alpha \neq \beta$. Then let $\phi$ be a substitution such that $x\mapsto [2]_{{\baxt}_2}$, we have
%\begin{align*}
%[1^{\alpha}2]_{\baxt_2}\cdot [u]_{\baxt_2}=\phi(\bfu[x])\neq \phi(\bfv[x])=[1^{\beta}2]_{\baxt_2}\cdot [v]_{\baxt_2}
%\end{align*}
%where $u,v \in \mathcal{A}_2^{\star}$, a contradiction. Thus $\alpha = \beta$. Therefore $\bfv[x]\in (x^*)^{\alpha}x\{x,x^*\}^{\times}$. Condition (Ib) holds by symmetry.

If $\bfu[x,y]\in x^{\alpha}x^*\{x,x^*,y,y^*\}^{\times}$,
then $\bfv[x,y]\in x^{\beta}z\{x,x^*,y,y^*\}^{\times}$ with $z\neq x$.
Suppose that $z\neq x^*$. Let $\phi_3$ be the substitution such that $x\mapsto [2]_{{\baxt}_2}, z\mapsto [1^{\alpha}2]_{{\baxt}_2}$. Then
\begin{align*}
[1^{\alpha}2]_{\baxt_2}\cdot [\bfs]_{\baxt_2}=\phi_3(\bfu[x,y])\neq \phi_3(\bfv[x,y])=[1^{\alpha+ \beta}2]_{\baxt_2}\cdot [\bft]_{\baxt_2}
\end{align*}
where $\bfs, \bft \in \mathcal{A}_2^{\star}$, a contradiction. Thus $z = x^*$. Suppose that $\alpha \neq \beta$. Let $\phi_4$ be the substitution such that $x\mapsto [2]_{{\baxt}_2}$. Then
\begin{align*}
[1^{\alpha}2]_{\baxt_2}\cdot [\bfs]_{\baxt_2}=\phi_4(\bfu[x,y])\neq \phi_4(\bfv[x,y])=[1^{\beta}2]_{\baxt_2}\cdot [\bft]_{\baxt_2}
\end{align*}
where $\bfs, \bft \in \mathcal{A}_2^{\star}$, a contradiction. Thus $\alpha = \beta$, and so $\bfv[x,y]\in x^{\alpha}x^*\{x,x^*,y,y^*\}^{\times}$. Therefore (Ia) holds.

If $\bfu[x,y]\in y^{\alpha}x\{x,x^*,y,y^*\}^{\times}$, then $\bfv[x,y]\in y^{\beta}z\{x,x^*,y,y^*\}^{\times}$ with $z\neq y$. Clearly $z\neq y^*$ by condition (Ia). If $z=x^*$, then $\bfu[x]$ starts with $x$ but $\bfv[x]$ starts with $x^*$, which is impossible.  Therefore $z=x$. Suppose that $\alpha \neq \beta$.
Let $\phi_5$ be the substitution such that $x\mapsto [2]_{{\baxt}_2}, y\mapsto [1]_{{\baxt}_2}$. Then
\begin{align*}
[1^{\alpha}2]_{\baxt_2}\cdot [\bfs]_{\baxt_2}=\phi_5(\bfu[x,y])\neq\phi_5(\bfv[x,y])=[1^{\beta}2]_{\baxt_2}\cdot [\bft]_{\baxt_2}
\end{align*}
where $\bfs, \bft \in \mathcal{A}_2^{\star}$, a contradiction. Thus $\alpha = \beta$, and so $\bfv[x,y]\in y^{\alpha}x\{x,x^*,y,y^*\}^{\times}$.
Therefore (IIa) holds.

If $\bfu[x,y]\in \bfa x\{x,x^*,y,y^*\}^{\times}$ with $\con(\bfa)= \{x^*, y^*\}$, then by condition (IIa) $\bfv[x,y]\in \bfa'y\{x,x^*,y,y^*\}^{\times}$ or
$\bfv[x,y]\in \bfa'x\{x,x^*,y,y^*\}^{\times}$ with $\con(\bfa')= \{x^*,y^*\}$. Suppose that $|\bfa'|\neq |\bfa|$. Let $\phi_6$ be the substitution from $\mathcal{X}$ to $\mathsf{baxt}_2$ such that $x\mapsto [2]_{{\baxt}_2}, y\mapsto [2]_{{\baxt}_2}$. Then
\begin{align*}
[1^{|\bfa|}2]_{\baxt_2}\cdot [\bfs]_{\baxt_2}=\phi_6(\bfu[x,y])\neq \phi_6(\bfv[x,y])=[1^{|\bfa'|}2]_{\baxt_2}\cdot [\bft]_{\baxt_2}
\end{align*}
where $\bfs, \bft \in \mathcal{A}_2^{\star}$, a contradiction. Therefore $|\bfa'|= |\bfa|$. Suppose that $\occ(x^*, \bfa')\neq \occ(x^*, \bfa)$. Let $\phi_7$ be the substitution from $\mathcal{X}$ to $\mathsf{baxt}_2$ such that $x\mapsto [22]_{{\baxt}_2}, y\mapsto [2]_{{\baxt}_2}$.  Then
 \begin{align*}
[1^{|\bfa|+\occ(x^*, \bfu)}2]_{\baxt_2}\cdot [\bfs]_{\baxt_2}=\phi_7(\bfu[x,y])\neq \phi_7(\bfv[x,y])=[1^{|\bfa'|+\occ(x^*, \bfv)}2]_{\baxt_2}\cdot [\bft]_{\baxt_2}
\end{align*}
where $\bfs, \bft \in \mathcal{A}_2^{\star}$, a contradiction.
Therefore $\occ(x^*, \bfa')= \occ(x^*, \bfa), \occ(y^*, \bfa')= \occ(y^*, \bfa)$, and so condition (IIIa) holds.

Conversely, let $\bfu\approx \bfv$ be any balanced word identity satisfying (I)--(III) and $\phi$ be any substitution from $\mathcal{X}$ to $\mathsf{baxt}_2$.
Note that $\ev(\phi(\bfu))=\ev(\phi(\bfv))$ since $\bfu \approx \bfv$ is balanced. If $\supp(\phi(\bfu))=\{1\}$ or $\{2\}$, then $\phi(\bfu)=\phi(\bfv)$.
If $\supp(\phi(\bfu))=\{1,2 \}$, to show $\phi(\bfu)=\phi(\bfv)$, it suffices to show $\rpi(\phi(\bfu))=\rpi(\phi(\bfv))$ and $\lpi(\phi(\bfu))=\lpi(\phi(\bfv))$. By symmetry, we only need to show that $\rpi(\phi(\bfu))=\rpi(\phi(\bfv))$.
If $\phi(\bfu)$ does not have a $2$-$1$ right precedence, then $\phi(\bfv)$ does not have a $2$-$1$ right precedence by conditions (I) and (II).
If $\phi(\bfu)$ has a $2$-$1$ right precedence of index $r$, then we may assume that $\phi(x)\neq [\varepsilon]_{\baxt_2}$ for any $x\in \con(\bfu)$ by Remark \ref{rem:delete}, whence $\bfu$ can be written into the form $\bfu = \bfu_1z\bfu_2$ for some possibly empty words $\bfu_1, \bfu_2$
satisfying $\supp(\phi(\bfu_2))=\{2\}$ and $\supp(\phi(z))=\{1, 2\}$ or
$\{1\}$.  Note that $z\not\in\con(\bfu_2)$ and $z^*\in\con(\bfu_2)$ if $\supp(\phi(z))=\{1\}$ and that for any $x \in \con(\bfu_2), x^* \not \in \con(\bfu_2)$. As such, $\phi(z\bfu_2)$ has the same $2$-$1$ right precedence of index $r$ as $\phi(\bfu)$.
Now we consider the form of $\bfv$. There are two cases.

{\bf Case~1.} $z^* \not \in \con(\bfu_2)$.  It follows from (IIb) that $\bfv = \bfv_1z\bfv_2$ where $\bfv_2$ is a permutation of $\bfu_2$.  This implies that $\supp(\phi(\bfv_2))=\{2\}$, hence $\phi(z\bfv_2)$ also has a $2$-$1$ right precedence of index $r$. Since $\phi(z\bfv_2)$ has the same right precedence as $\phi(\bfv)$, $\phi(\bfv)$ has a $2$-$1$ right precedence of index $r$. Therefore $\rpi(\phi(\bfu))=\rpi(\phi(\bfv))$.

{\bf Case~2.} $z^* \in \con(\bfu_2)$. Let $\con(\bfu_2)=\{z^*, x_1, x_2, \dots, x_n\}$.  Without loss of generality, we may assume that $\fp(\bfu_2)=z^* x_1x_2\cdots x_n$. Note that $x_1^*, x_2^*, \dots, x_n^*\not\in \con(\bfu_2)$.  Then $\bfv=\bfv_1\bfy z^*\bfx_1x_1\bfx_2x_2\cdots\bfx_nx_n$ by (IIb) where $\bfy\in \{x_1,x_1^*, \dots, x_n, x_n^*, z, z^*\}^{+}$ and $\bfx_i \in \{x_i, x_{i+1},\dots, x_n\}^{+}$ for $i=1,2,\dots, n$.
%Further, we have $\bfx_i \in \{x_i, x_{i+1},$ $\dots, x_n\}^{+}$ for $i=1,2,\dots, n$.
%Otherwise, suppose that $x_j^* \in \con(\bfx_i)$ with $i\leq j\leq n$. Then $\bfu[x_{i-1}, x_j]\in \{x_{i-1}, x_{i-1}^*, x_j, x_j^*\}^{+}x_{i-1}x_j^{\alpha}$ but $\bfv[x_{i-1}, x_j]\in \{x_{i-1}, x_{i-1}^*, x_j, x_j^*\}^{+}$ $x_{j}^{*}x_j^{\alpha}$ when $i\geq 2$ and $\bfu[z, x_j]\in \{z, z^*, x_j, x_j^*\}^{+}z^*x_j^{\alpha}$ but $\bfv[z, x_j]\in \{z, z^*,  x_j, x_j^*\}^{+}x_{j}^{*}x_j^{\alpha}$ when $i=1$, which contradicts conditions (IIb) and (IIIb).
Therefore either $\bfv=\bfv_1z\bfv_2$ or $\bfv=\bfv_1x_i^*\bfv_2$ with $\con(\bfv_2)= \con(\bfu_2)$.

{\bf 2.1.} $\bfv=\bfv_1z\bfv_2$. It follows from (IIIb) that $\bfv_2$ is a permutation of $\bfu_2$. This implies that $\supp(\phi(\bfv_2))=\{2\}$, hence $\phi(z\bfv_2)$ also has a $2$-$1$ right precedence of index $r$.
Since $\phi(z\bfv_2)$ has the same right precedence as $\phi(\bfv)$, $\phi(\bfv)$ has a $2$-$1$ right precedence of index $r$.
Therefore $\rpi(\phi(\bfu))=\rpi(\phi(\bfv))$.

{\bf 2.2.} $\bfv=\bfv_1x_i^*\bfv_2$. It follows from (IIIb) that $\occ(z^*, \bfu_2)=\occ(z^*, \bfv_2)$ and  $\occ(x_i, \bfu_2)=\occ(x_i, \bfv_2)$. Suppose that there exists some $x_j\in\con(\bfu_2)$ such that $\occ(x_j, \bfu_2)\neq\occ(x_j, \bfv_2)$. Then $x_j\neq x_i, z^*$.
If $x_j^* \not\in \con(\bfu_1)$, then it follows from (IIIb) that $\overrightarrow{\occ}_{x_i^*}(x_j,\bfu)=\overrightarrow{\occ}_{x_i^*}(x_j,\bfv)$ and $\overrightarrow{\occ}_{z}(x_j,\bfu)=\overrightarrow{\occ}_{z}(x_j,\bfv)$.  It follows from the forms of $\bfu$ and $\bfv$ that $\overrightarrow{\occ}_{z}(x_j,\bfv)-\overrightarrow{\occ}_{x_i^*}(x_j,\bfv)\geq 0$ and  $\overrightarrow{\occ}_{z}(x_j,\bfu)-\overrightarrow{\occ}_{x_i^*}(x_j,\bfu)\leq0$. Suppose  that $\overrightarrow{\occ}_{z}(x_j,\bfv)-\overrightarrow{\occ}_{x_i^*}(x_j,\bfv)>0$. Then $\overrightarrow{\occ}_{z}(x_j,\bfu)-\overrightarrow{\occ}_{x_i^*}(x_j,\bfu)>0$, which is impossible.
If $x_j^* \in \con(\bfu_1)$, then there are four cases about the order of $_{\infty}x_j^*, {_{\infty}x_i^*}, {_{\infty}z}$ in $\bfu, \bfv$:
$_{\infty}x_j^* \prec_{\bfu} {_{\infty}x_i^*} \prec_{\bfu} {_{\infty}z}$ and $_{\infty}x_j^* \prec_{\bfv} {_{\infty}z} \prec_{\bfv} {_{\infty}x_i^*}$;
$_{\infty}x_j^* \prec_{\bfu} {_{\infty}x_i^*} \prec_{\bfu} {_{\infty}z}$ and $  {_{\infty}z} \prec_{\bfv} {_{\infty}x_j^*} \prec_{\bfv} {_{\infty}x_i^*}$;
$ {_{\infty}x_i^*}\prec_{\bfu} {_{\infty}x_j^*} \prec_{\bfu} {_{\infty}z}$ and $_{\infty}x_j^* \prec_{\bfv} {_{\infty}z} \prec_{\bfv} {_{\infty}x_i^*}$;
and $ {_{\infty}x_i^*}  \prec_{\bfu} {_{\infty}x_j^*} \prec_{\bfu} {_{\infty}z}$ and $  {_{\infty}z} \prec_{\bfv} {_{\infty}x_j^*} \prec_{\bfv}$.

If $_{\infty}x_j^* \prec_{\bfu} {_{\infty}x_i^*} \prec_{\bfu} {_{\infty}z}$ and $_{\infty}x_j^* \prec_{\bfv} {_{\infty}z} \prec_{\bfv} {_{\infty}x_i^*}$, then it follows from (IIIb) that $\overrightarrow{\occ}_{x_i^*}(x_j,\bfu)=\overrightarrow{\occ}_{x_i^*}(x_j,\bfv)$ and $\overrightarrow{\occ}_{z}(x_j,\bfu)=\overrightarrow{\occ}_{z}(x_j,\bfv)$. By the forms of $\bfu$ and $\bfv$,
$\overrightarrow{\occ}_{z}(x_j,\bfv)-\overrightarrow{\occ}_{x_i^*}(x_j,\bfv)\geq0$, $\overrightarrow{\occ}_{z}(x_j,\bfu)-\overrightarrow{\occ}_{x_i^*}(x_j,\bfu)\leq0$.
Suppose $\overrightarrow{\occ}_{z}(x_j,\bfv)-\overrightarrow{\occ}_{x_i^*}(x_j,\bfv)>0$.  Then $\overrightarrow{\occ}_{z}(x_j,\bfu)-\overrightarrow{\occ}_{x_i^*}(x_j,\bfu)>0$, which is impossible.

If $_{\infty}x_j^* \prec_{\bfu} {_{\infty}x_i^*} \prec_{\bfu} {_{\infty}z}$ and $  {_{\infty}z} \prec_{\bfv} {_{\infty}x_j^*} \prec_{\bfv} {_{\infty}x_i^*}$, then it follows from (IIIb) that $\overrightarrow{\occ}_{x_i^*}(x_j,\bfu)=\overrightarrow{\occ}_{x_i^*}(x_j,\bfv)$ and $\overrightarrow{\occ}_{z}(x_j,\bfu)=\overrightarrow{\occ}_{x_j^*}(x_j,\bfv)$. By the forms of $\bfu$ and $\bfv$, $\overrightarrow{\occ}_{x_j^*}(x_j,\bfv)-\overrightarrow{\occ}_{x_i^*}(x_j,\bfv)\geq0$, $\overrightarrow{\occ}_{z}(x_j,\bfu)-\overrightarrow{\occ}_{x_i^*}(x_j,\bfu)\leq0$. Suppose  $\overrightarrow{\occ}_{x_j^*}(x_j,\bfv)-\overrightarrow{\occ}_{x_i^*}(x_j,\bfv)>0$. Then $\overrightarrow{\occ}_{z}(x_j,\bfu)-\overrightarrow{\occ}_{x_i^*}(x_j,\bfu)>0$, which is impossible.

If $ {_{\infty}x_i^*}\prec_{\bfu} {_{\infty}x_j^*} \prec_{\bfu} {_{\infty}z}$ and $_{\infty}x_j^* \prec_{\bfv} {_{\infty}z} \prec_{\bfv} {_{\infty}x_i^*}$, then it follows from (IIIb) that $\overrightarrow{\occ}_{x_j^*}(x_j,\bfu)=\overrightarrow{\occ}_{x_i^*}(x_j,\bfv)$ and  $\overrightarrow{\occ}_{z}(x_j,\bfu)=\overrightarrow{\occ}_{z}(x_j,\bfv)$. By the forms of $\bfu$ and $\bfv$, $\overrightarrow{\occ}_{z}(x_j,\bfv)-\overrightarrow{\occ}_{x_i^*}(x_j,\bfv)\geq0$, $\overrightarrow{\occ}_{z}(x_j,\bfu)-\overrightarrow{\occ}_{x_j^*}(x_j,\bfu)\leq0$. Suppose $\overrightarrow{\occ}_{z}(x_j,\bfv)-\overrightarrow{\occ}_{x_i^*}(x_j,\bfv)>0$. Then $\overrightarrow{\occ}_{z}(x_j,\bfu)-\overrightarrow{\occ}_{x_j^*}(x_j,\bfu)>0$, which is impossible.

If $ {_{\infty}x_i^*}  \prec_{\bfu} {_{\infty}x_j^*} \prec_{\bfu} {_{\infty}z}$ and $  {_{\infty}z} \prec_{\bfv} {_{\infty}x_j^*} \prec_{\bfv} {_{\infty}x_i^*}$, then it follows from (IIIb) that $\overrightarrow{\occ}_{x_j^*}(x_j,\bfu)=\overrightarrow{\occ}_{x_i^*}(x_j,\bfv)$ and $\overrightarrow{\occ}_{z}(x_j,\bfu)=\overrightarrow{\occ}_{x_j^*}(x_j,\bfv)$. By the forms of $\bfu$ and $\bfv$, $\overrightarrow{\occ}_{x_j^*}(x_j,\bfv)-\overrightarrow{\occ}_{x_i^*}(x_j,\bfv)\geq0$, $\overrightarrow{\occ}_{z}(x_j,\bfu)-\overrightarrow{\occ}_{x_j^*}(x_j,\bfu)\leq0$.  Suppose $\overrightarrow{\occ}_{x_j^*}(x_j,\bfv)-\overrightarrow{\occ}_{x_i^*}(x_j,\bfv)>0$. Then $\overrightarrow{\occ}_{z}(x_j,\bfu)-\overrightarrow{\occ}_{x_j^*}(x_j,\bfu)>0$, which is impossible.

Thus $\bfv = \bfv_1x_i^*\bfv_2$ where $\bfv_2$ is a permutation of $\bfu_2$. This implies that $\phi(\bfv_2)$ has support $\{2\}$, hence $\phi(x_i^*\bfv_2)$ also has a $2$-$1$ right precedence of index $r$.
Since $\phi(x_i^*\bfv_2)$ has the same right precedence as $\phi(\bfv)$, $\phi(\bfv)$ has a $2$-$1$ right precedence of index $r$. Therefore $\rpi(\phi(\bfu))=\rpi(\phi(\bfv))$.
\end{proof}

\begin{theorem}\label{thm:baxt3}
A word identity $\bfu\approx \bfv$ holds in $(\mathsf{baxt}_3,~^\sharp)$ if and only if $\bfu\approx \bfv$ is balanced and satisfies that for any $x, y \in \con(\bfu)$ with $x, x^*\neq y$,
\begin{enumerate}[\rm(I)]
\item \begin{enumerate}[\rm(a)]
          \item  if $\bfu[x,y]\in x^{\alpha}x^*\{x,x^*,y,y^*\}^{\times}$ for some $\alpha\geq 1$, then $\bfv[x,y]\in x^{\alpha}x^*\{x,x^*,y,y^*\}^{\times}$,

          \item if $\bfu[x,y]\in \{x,x^*,y,y^*\}^{\times}x^*x^{\alpha}$ for some $\alpha\geq 1$, then $\bfv[x,y]\in \{x,x^*,y,y^*\}^{\times}x^*x^{\alpha}$.
        \end{enumerate}

  \item \begin{enumerate}[\rm(a)]
          \item  if $\bfu[x,y]\in y^{\alpha}x\{x,x^*,y,y^*\}^{\times}$ for some $\alpha\geq 1$, then $\bfv[x,y]\in y^{\alpha}x\{x,x^*,y,y^*\}^{\times}$,

          \item if $\bfu[x,y]\in \{x,x^*,y,y^*\}^{\times}xy^{\alpha}$ for some $\alpha\geq 1$, then $\bfv[x,y]\in \{x,x^*,y,y^*\}^{\times}xy^{\alpha}$.
        \end{enumerate}

\item \begin{enumerate}[\rm(a)]
          \item  if $\bfu[x,y]\in \bfa x\{x,x^*,y,y^*\}^{\times}$ with $\con(\bfa)= \{x^*, y^*\}$, then
          $\bfv[x,y]\in \bfa' x\{x,x^*,y,y^*\}^{\times}$
            where $\bfa'$ is a permutation of $\bfa$,
          \item   if $\bfu[x,y]\in \{x, x^*, y, y^*\}^{\times} x\bfa$ with $\con(\bfa)= \{x^*, y^*\}$, then
          $\bfv[x,y]\in \{x,x^*,y,y^*\}^{\times}x\bfa'$
             where $\bfa'$ is a permutation of $\bfa$.
        \end{enumerate}

\item \begin{enumerate}[\rm(a)]
          \item $\overleftarrow{\occ}_y(x, \bfu)+\overleftarrow{\occ}_y(x^*, \bfu)= \overleftarrow{\occ}_y(x, \bfv)+\overleftarrow{\occ}_y(x^*, \bfv)$,
           \item $\overrightarrow{\occ}_y(x, \bfu)+\overrightarrow{\occ}_y(x^*, \bfu)= \overrightarrow{\occ}_y(x, \bfv)+\overrightarrow{\occ}_y(x^*, \bfv)$.
        \end{enumerate}
 \item \begin{enumerate}[\rm(a)]
          \item  if $\bfu[x,y]\in \bfa y\{x,x^*,y,y^*\}^{\times}$ with $\con(\bfa)= \{x, x^*\}$, then
          $\bfv[x,y]\in \bfa' y\{x,x^*,y,y^*\}^{\times}$
           where $\bfa'$ is a permutation of $\bfa$,
          \item   if $\bfu[x,y]\in \{x, x^*, y, y^*\}^{\times} y\bfa$ with $\con(\bfa)= \{x, x^*\}$, then
          $\bfv[x,y]\in \{x,x^*,y,y^*\}^{\times}y\bfa'$
            where $\bfa'$ is a permutation of $\bfa$.
        \end{enumerate}
\end{enumerate}
\end{theorem}

%\begin{remark}
%It is easy to see that
%\begin{enumerate}[\rm(a)]
%  \item Condition (IVa) is equivalent to that the longest prefix of $\bfu[x,y]$ in which $x$ can not occur simultaneously and $y, y^*$ can not occur simultaneously has the same content and number of occurrences as the longest prefix of $\bfv[x,y]$ in which $x$ can not occur simultaneously and $y, y^*$ can not occur simultaneously and that Condition (IVb) is equivalent to that the longest suffix of $\bfu[x,y]$ in which $x$ can not occur simultaneously and $y, y^*$ can not occur simultaneously has the same content and number of occurrences longest suffix of $\bfv[x,y]$ in which $x$ can not occur simultaneously and $y, y^*$ can not occur simultaneously.
%\end{enumerate}
%\end{remark}

\begin{proof}
%First, we prove that a word identity $\bfu\approx \bfv$ satisfied by $(\mathsf{baxt}_3,~^\sharp)$ must be  balanced and satisfy (I)--(V).
Suppose that $\bfu \approx  \bfv$ is a word identity satisfied by $(\mathsf{baxt}_3,~^\sharp)$. It is routine to verify that the involution submonoid of $(\mathsf{baxt}_3,~^\sharp)$ generated by $[1]_{\mathsf{baxt}_3}$ and $[3]_{\mathsf{baxt}_3}$  is isomorphic to $(\mathsf{baxt}_2,~^\sharp)$. It follows from Theorem~\ref{thm:baxt2} that $\bfu\approx \bfv$ is balanced and satisfies conditions (I)--(II). Now to show that (IIIa)---(Va) hold, and then (IIIb)---(Vb) hold by symmetry.

Suppose that $\bfu[x,y]\in \bfa x\{x,x^*,y,y^*\}^{\times}$ with $\con(\bfa)= \{x^*, y^*\}$. Then it follows from condition (IIIa) of Theorem~\ref{thm:baxt2} that
$\bfv[x,y]\in \bfa' x\{x,x^*,y,y^*\}^{\times}$ or $\bfv[x,y]\in \bfa' y\{x,x^*,y,y^*\}^{\times}$ where $\bfa'$ is a permutation of $\bfa$. Suppose that $\bfv[x,y]\in \bfa' y\{x,x^*,y,y^*\}^{\times}$. Let $\phi_1$ be a substitution such that $x\mapsto [3]_{{\baxt}_3}, y\mapsto [2]_{{\baxt}_3}$ and any other variable to $[\varepsilon]_{{\baxt}_3}$. Then $\phi_1(\bfu)$ have a $2$-$3$ left precedence of index $\occ(y^*, \bfa)$ but the index of $2$-$3$ left precedence of $\phi_1(\bfv)$ is greater than $\occ(y^*, \bfa)$, a contradiction. Therefore $\bfv[x,y]\in \bfa' x\{x,x^*,y,y^*\}^{\times}$, and so (IIIa) holds.

Suppose that $\overleftarrow{\occ}_y(x, \bfu)+\overleftarrow{\occ}_y(x^*, \bfu)\neq \overleftarrow{\occ}_y(x, \bfv)+\overleftarrow{\occ}_y(x^*, \bfv)$ for some $x,x^*\neq y \in \con(\bfu)$. Let $\phi_2$ be a substitution such that $x\mapsto [2]_{{\baxt}_3}, y\mapsto [3]_{{\baxt}_3}$ and any other variable to $[\varepsilon]_{{\baxt}_3}$. Then $\phi_2(\bfu)$ has a $2$-$3$ left precedence of index $\overleftarrow{\occ}_y(x, \bfu)+\overleftarrow{\occ}_y(x^*, \bfu)$ but $\phi_2(\bfv)$ has a $2$-$3$ left precedence of index $\overleftarrow{\occ}_y(x, \bfv)+\overleftarrow{\occ}_y(x^*, \bfv)$, a contradiction. Therefore $\overleftarrow{\occ}_y(x, \bfu)+\overleftarrow{\occ}_y(x^*, \bfu)= \overleftarrow{\occ}_y(x, \bfv)+\overleftarrow{\occ}_y(x^*, \bfv)$, and so (IVa) holds.

If $\bfu[x,y]\in \bfa y\{x,x^*,y,y^*\}^{\times}$ with $\con(\bfa)= \{x, x^*\}$, then it follows from (I) that $\bfv[x,y]\in \bfa' y\{x,x^*,y,y^*\}^{\times}$ with $\con(\bfa')= \{x, x^*\}$. By (IVb), $\occ(x, \bfa)+\occ(x^*, \bfa)=\occ(x, \bfa')+\occ(x^*, \bfa')$.
Suppose that $\occ(x, \bfa)\neq\occ(x, \bfa')$. Let $\phi_3$ be a substitution such that $x\mapsto [1]_{{\baxt}_3}, y\mapsto [2]_{{\baxt}_3}$ and any other variable to $[\varepsilon]_{{\baxt}_3}$. Then $\phi_3(\bfu)$ has a $1$-$2$ left precedence of index $\occ(x, \bfa)$ but $\phi_3(\bfv)$ has a $1$-$2$ left precedence of index $\occ(x, \bfa')$, a contradiction. Therefore $\bfv[x,y]\in \bfa' y\{x,x^*,y,y^*\}^{\times}$
where $\bfa'$ is a permutation of $\bfa$, and so (Va) holds.

Conversely, let $\bfu\approx \bfv$ be any balanced word identity satisfying (I)--(V) and $\phi$ be any substitution from $\mathcal{X}$ to $\mathsf{baxt}_3$.
Note that $\ev(\phi(\bfu))=\ev(\phi(\bfv))$ since $\bfu \approx \bfv$ is balanced. Now, to show $\phi(\bfu)=\phi(\bfv)$, it suffices to show $\rpi(\phi(\bfu))=\rpi(\phi(\bfv))$ and $\lpi(\phi(\bfu))=\lpi(\phi(\bfv))$. By symmetry, we only need to show that $\rpi(\phi(\bfu))=\rpi(\phi(\bfv))$. If $\phi(\bfu)$ does not have a $2$-$1$ [resp. $3$-$2$, $3$-$1$] right precedence, then $\phi(\bfv)$ does not have a $2$-$1$ [resp. $3$-$2$, $3$-$1$] right precedence by (I) and (II). In the following, we show that if $\phi(\bfu)$ has a $2$-$1$ [resp. $3$-$2$, $3$-$1$] right precedence of index $r$, then $\phi(\bfv)$ has a $2$-$1$ [resp. $3$-$2$, $3$-$1$] right precedence of index $r$.

If $\phi(\bfu)$ has a $3$-$1$ right precedence of index $r$, then we may assume that $\phi(x)\neq [\varepsilon]_{\baxt_3}$ for any $x\in \con(\bfu)$ by Remark \ref{rem:delete}, whence $\bfu$ can be written into the form $\bfu = \bfu_1z\bfu_2$ for some possibly empty words $\bfu_1, \bfu_2$
satisfying $\supp(\phi(\bfu_2))=\{3\}$ and $\supp(\phi(z))=\{1\}$ or $\{1, 3\}$. Note that $z\not\in\con(\bfu_2)$ and $z^*$ may occur in $\bfu_2$ if $\supp(\phi(z))=\{1\}$ and that for any $x \in \con(\bfu_2), x^* \not \in \con(\bfu_2)$. As such, $\phi(z\bfu_2)$ has the same $3$-$1$ right precedence of index $r$ as $\phi(\bfu)$.
 Now we consider the form of $\bfv$.
If $z^* \not \in \con(\bfu_2)$, then it follows from (IIb) that $\bfv = \bfv_1z\bfv_2$ where $\bfv_2$ is a permutation of $\bfu_2$. This implies that $\supp(\phi(\bfv_2))=\{3\}$, hence $\phi(z\bfv_2)$ also has a $3$-$1$ right precedence of index $r$. Since $\phi(z\bfv_2)$ has the same $3$-$1$ right precedence of index $r$ as $\phi(\bfv)$, $\phi(\bfv)$ has a $3$-$1$ right precedence of index $r$.
If $z^* \in \con(\bfu_2)$. Let $\con(\bfu_2)=\{z^*, x_1, x_2, \dots, x_n\}$. Without loss of generality, we may assume that $\fp(\bfu_2)=z^* x_1x_2\cdots x_n$. Note that $x_1^*, x_2^*, \dots, x_n^*\not\in \con(\bfu_2)$.  Then it follows from (IIb) that $\bfv=\bfv_1\bfy z^*\bfx_1x_1$ $\bfx_2x_2\cdots\bfx_nx_n$ where $\bfy\in \{x_1,x_1^*, \dots, x_n, x_n^*, z, z^*\}^{+}$ and $\bfx_i \in \{x_i, x_{i+1}, \dots, x_n\}^{+}$ for $i=1,2,\dots, n$.
%Further, we have $\bfx_i \in \{x_i, x_{i+1},$ $ \dots, x_n\}^{+}$ for $i=1,2,\dots, n$.
%Otherwise, suppose that $x_j^* \in \con(\bfx_i)$ with $i\leq j\leq n$. Then $\bfu[x_{i-1}, x_j]\in \{x_{i-1}, x_{i-1}^*, x_j, x_j^*\}^{+}x_{i-1}x_j^{\alpha}$ but $\bfv[x_{i-1}, x_j]\in \{x_{i-1}, x_{i-1}^*, x_j, x_j^*\}^{+}$ $x_{j}^{*}x_j^{\alpha}$ when $i\geq 2$ and $\bfu[z, x_j]\in \{z, z^*, x_j, x_j^*\}^{+}z^*x_j^{\alpha}$ but $\bfv[z, x_j]\in \{z, z^*,  x_j, x_j^*\}^{+}x_{j}^{*}x_j^{\alpha}$ when $i=1$, which contradicts conditions (IIb) and (IIIb).
Therefore $\bfv=\bfv_1z\bfv_2$ where $\bfv_2$ is a permutation of $\bfu_2$ by (IIIb). This implies that $\phi(\bfv_2)$ has support $\{3\}$, hence $\phi(z\bfv_2)$ also has a $3$-$1$ right precedence of index $r$.
Since $\phi(z\bfv_2)$ has the $3$-$1$ right precedence of index $r$ as $\phi(\bfv)$, $\phi(\bfv)$ has a $3$-$1$ right precedence of index $r$.

If $\phi(\bfu)$ has a $2$-$1$ right precedence of index $r$, then we may assume that $\phi(x)\neq [\varepsilon]_{\baxt_3}$ for any $x\in \con(\bfu)$ by Remark \ref{rem:delete}, whence $\bfu$ can be written into the form $\bfu = \bfu_1z\bfu_2$ for
some possibly empty words $\bfu_1, \bfu_2$
satisfying $\supp(\phi(\bfu_2))\subseteq\{2, 3\}$ and $\supp(\phi(z))=\{1\}$,
$\{1, 2\}$, $\{1, 3\}$ or $\{1, 2, 3\}$. Note that $z\not\in\con(\bfu_2)$ and $z^*$ may occur in $\bfu_2$ if $\supp(\phi(z))=\{1\}$ or $\{1, 2\}$ and that for any $x \in \con(\bfu_2)$, $x^*$ may occur in $\con(\bfu_2)$ if $\supp(\phi(x))=\{2\}$.  As such, $\phi(z\bfu_2)$ has the same $2$-$1$ right precedence of index $r$ as $\phi(\bfu)$.
It follows from (IIb)--(Vb) that $\bfv=\bfv_1z\bfv_2$ satisfying that $z\not \in \con(\bfv_2)$, $\con(\bfu_2)=\con(\bfv_2)$, and $\occ(z^*, \bfu_2)=\occ(z^*, \bfv_2)$, and $\occ(x, \bfu_2)=\occ(x, \bfv_2)$ when $x\in \con(\bfu_2), x^*\not\in \con(\bfu_2)$, and $\occ(x, \bfu_2)+\occ(x^*, \bfu_2)= \occ(x, \bfv_2)+\occ(x^*, \bfv_2)$ when $x, x^*\in \con(\bfu_2)$.
This implies that $\phi(z\bfv_2)$ also has a $2$-$1$ right precedence of index $r$.
Since $\phi(z\bfv_2)$ has the same $2$-$1$ right precedence of index $r$ as $\phi(\bfv)$, $\phi(\bfv)$ has a $2$-$1$ right precedence of index $r$.

If $\phi(\bfu)$ has a $3$-$2$ right precedence of index $r$, then we may assume that $\phi(x)\neq [\varepsilon]_{\baxt_3}$ for any $x\in \con(\bfu)$ by Remark \ref{rem:delete}, whence $\bfu$ can be written into the form $\bfu = \bfu_1z\bfu_2$
satisfying $\supp(\phi(\bfu_2))\subseteq\{1, 3\}$ and $\supp(\phi(z))=\{2\}$,
$\{1, 2\}$, $\{2, 3\}$ or $\{1, 2, 3\}$. Note that $z, z^*$ cannot occur in $\bfu_2$ and that for any $x \in \con(\bfu_2)$, $x^*$ may occur in $\con(\bfu_2)$. As such, $\phi(z\bfu_2)$ has the same $3$-$2$ right precedence of index $r$ as $\phi(\bfu)$.
It follows from (IIb)-(Vb) that $\bfv=\bfv_1z\bfv_2$ such that $z\not \in \con(\bfv_2)$, $\con(\bfu_2)=\con(\bfv_2)$, $\occ(x, \bfu_2)=\occ(x, \bfv_2)$ and $\occ(x^*, \bfu_2)= \occ(x^*, \bfv_2)$.
This implies that $\phi(z\bfv_2)$ also has a $3$-$2$ right precedence of index $r$.
Since $\phi(z\bfv_2)$ has the same $3$-$2$ right precedence of index $r$ as $\phi(\bfv)$, $\phi(\bfv)$ has a $3$-$2$ right precedence of index $r$.
Therefore $\rpi(\phi(\bfu))=\rpi(\phi(\bfv))$.
\end{proof}

\begin{theorem}\label{thm:baxt4+}
A word identity $\bfu\approx \bfv$ holds in $(\mathsf{baxt}_4,~^\sharp)$ if and only if $\bfu\approx \bfv$ is balanced and $\overrightarrow{\occ}_x(y, \bfu)=\overrightarrow{\occ}_x(y, \bfv), \overleftarrow{\occ}_x(y, \bfu)=\overleftarrow{\occ}_x(y, \bfv)$ for any $x, y \in \con(\bfu)$.
\end{theorem}

\begin{proof}
%First, we prove that a word identity $\bfu\approx \bfv$ satisfied by $(\mathsf{baxt}_4,~^\sharp)$ must be  balanced and $\overrightarrow{\occ}_x(y, \bfu)=\overrightarrow{\occ}_x(y, \bfv), \overleftarrow{\occ}_x(y, \bfu)=\overleftarrow{\occ}_x(y, \bfv)$.
Suppose that $\bfu \approx  \bfv$ is a word identity satisfied by $(\mathsf{baxt}_4,~^\sharp)$. It is routine to verify that the involution submonoid of $(\mathsf{baxt}_4,~^\sharp)$ generated by $[1]_{\mathsf{baxt}_4}$ and $[4]_{\mathsf{baxt}_3}$  is isomorphic to $(\mathsf{baxt}_2,~^\sharp)$. Then it follows from Theorem~\ref{thm:baxt2} that $\bfu\approx \bfv$ is balanced.
Suppose that $\overrightarrow{\occ}_x(y, \bfu)\neq\overrightarrow{\occ}_x(y, \bfv)$. Let $\phi$ be a substitution such that $x\mapsto [1]_{{\baxt}_4}, y\mapsto [2]_{{\baxt}_4}$ and any other variable to $[\varepsilon]_{{\baxt}_4}$. Then $\phi(\bfu)$ has a 2-1 right precedence of index $\overrightarrow{\occ}_x(y, \bfu)$ but $\phi(\bfv)$ has a 2-1 left precedence of index $\overrightarrow{\occ}_x(y, \bfv)$, a contradiction. Therefore $\overrightarrow{\occ}_x(y, \bfu)= \overrightarrow{\occ}_x(y, \bfv)$ for any $x, y \in \con(\bfu)$.
Symmetrically, $\overleftarrow{\occ}_x(y, \bfu)=\overleftarrow{\occ}_x(y, \bfv)$ for any $x, y \in \con(\bfu)$.

Conversely, let $\bfu\approx \bfv$ be any balanced word identity  satisfying $\overrightarrow{\occ}_x(y, \bfu)=\overrightarrow{\occ}_x(y, \bfv), \overleftarrow{\occ}_x(y, \bfu)=\overleftarrow{\occ}_x(y, \bfv)$ for any $x, y \in \con(\bfu)$ and $\phi$ be any substitution from $\mathcal{X}$ to $\mathsf{baxt}_4$.
Note that $\ev(\phi(\bfu))=\ev(\phi(\bfv))$ since $\bfu \approx \bfv$ is balanced. Now, to show $\phi(\bfu)=\phi(\bfv)$, it suffices to show that $\rpi(\phi(\bfu))=\rpi(\phi(\bfv))$, and then $\lpi(\phi(\bfu))=\lpi(\phi(\bfv))$ by symmetry. If $\phi(\bfu)$ does not have a $j$-$i$ right precedence with $i<j$, then $\phi(\bfv)$ does not have a $j$-$i$ right precedence by $\overrightarrow{\occ}_x(y, \bfu)=\overrightarrow{\occ}_x(y, \bfv)$.
If $\phi(\bfu)$ has a $j$-$i$ right precedence of index $r$, then $\bfu$ can be written into the form $\bfu = \bfu_1z\bfu_2$
satisfying $i \not\in \supp(\phi(\bfu_2))$ and $i \in \supp(\phi(z))$. Note that $z$ cannot occur in $\bfu_2$. As such, $\phi(z\bfu_2)$ has the same $j$-$i$ right precedence of index $r$ as $\phi(\bfu)$.
It follows from $\overrightarrow{\occ}_x(y, \bfu)=\overrightarrow{\occ}_x(y, \bfv)$ that $\bfv=\bfv_1z\bfv_2$ such that $z\not \in \con(\bfv_2)$ and $\bfv_2$ is a permutation of $\bfu_2$.
This implies that $\phi(z\bfv_2)$ also has a $j$-$i$ right precedence of index $r$.
Since $\phi(z\bfv_2)$ has the same $j$-$i$ right precedence of index $r$ as $\phi(\bfv)$, $\phi(\bfv)$ has a $j$-$i$ right precedence of index $r$.
Therefore $\rpi(\phi(\bfu_2))=\rpi(\phi(\bfv))$.
\end{proof}

\begin{remark}
It follows from Theorems~\ref{thm:baxt2}--\ref{thm:baxt4+} that the involution monoids $(\baxt_2,~^{\sharp})$, $(\baxt_3,~^{\sharp})$ and $(\baxt_n,~^{\sharp})$ with $n\geq 4$ generate different varieties. However, the monoids $\baxt_n$ with $n\geq 2$ generate the same variety by \cite[Theorem~3.12]{CMR21b} or \cite[Proposition~6.10]{CJKM21}. And by Theorems~\ref{thm:baxt2}--\ref{thm:baxt4+}, an identity $\bfu\approx \bfv$ holds in $\baxt_n$ with $n\geq 2$ if and only if
\begin{enumerate}[\rm(i)]
  \item $\bfu\approx \bfv$ is balanced;
  \item $\overrightarrow{\occ}_x(y, \bfu)=\overrightarrow{\occ}_x(y, \bfv)$ and  $\overleftarrow{\occ}_x(y, \bfu)=\overleftarrow{\occ}_x(y, \bfv)$ for any $x, y \in \con(\bfu)$.
\end{enumerate}
This fact also has been shown in \cite[Theorem~4.3]{CMR21b}.
\end{remark}

\section{Finite basis problem for $(\mathsf{baxt}_n,~^\sharp)$} \label{sec:FBP4} %
In this section, the finite basis problem for involution monoid $(\mathsf{baxt}_n,~^\sharp)$ for each finite $n$ is solved.
Clearly, $(\mathsf{baxt}_1,~^\sharp)$ is finitely based since its involution is trivial and it is commutative.

For any word $\bfu \in(\mathcal{A}\cup \mathcal{A}^{\star})^{+}$, denote by $\overleftarrow{\bfu}$  the word obtained from $\bfu$ by writing $\bfu$ in reverse. For example, if $\bfu=x^5zy^*(z^*)^3x^*$, then $\overleftarrow{\bfu}=x^*(z^*)^3y^*zx^5$. The identity $\bf\overleftarrow{\bfu}\approx \overleftarrow{\bfv}$ is called the \textit{reverse} of $\bfu\approx \bfv$.  It follows from the proof of \cite[Lemma~7(i)]{Lee17a} that if an involution semigroup $(S,\op)$ satisfies the word identity $\bfu\approx \bfv$, then $(S,\op)$ also satisfies $\bf\overleftarrow{\bfu}\approx \overleftarrow{\bfv}$.

First, we show that $(\mathsf{baxt}_2,~^\sharp)$ is finitely based.

\begin{theorem}\label{thm:baxt2's basis}
The identities \eqref{id: inv} and
\begin{align}
\left. \begin{array}[c]{rcl}
x^*hxkxysx^*tx  \approx x^*hxkyxsx^*tx, & x^*hxkxysxtx^*  \approx x^*hxkyxsxtx^*,  \\
 xhx^*kxysx^*tx  \approx xhx^*kyxsx^*tx, & xhx^*kxysxtx^*  \approx xhx^*kyxsxtx^*,
\end{array} \quad\right\} \label{e1} \\%
\left. \begin{array}[c]{rcl}
x^*hxkxysy^*ty  \approx x^*hxkyxsy^*ty,  & x^*hxkxysyty^*  \approx x^*hxkyxsyty^*,\\
xhx^*kxysy^*ty  \approx xhx^*kyxsy^*ty,  & xhx^*kxysyty^*  \approx xhx^*kyxsyty^*,
\end{array} \quad \right\} \\ %
xhykxysxty  \approx  xhykyxsxty, \quad  xhykxysytx  \approx xhykyxsytx,\qquad\qquad\\
xhykxysx^*ty^* \approx  xhykyxsx^*ty^*, \quad  xhykxysy^*tx^*  \approx xhykyxsy^*tx^*,\qquad\\
x^*hy^*kxysx^*ty^*  \approx   x^*hy^*kyxsx^*ty^*, \,\,  x^*hy^*kxysy^*tx^*  \approx   x^*hy^*kyxsy^*tx^*,\\
\left. \begin{array}[c]{rcl}
x^*hxkxysxty  \approx x^*hxkyxsxty,  & x^*hxkxysytx  \approx x^*hxkyxsytx, \\
xhx^*kxysxty  \approx xhx^*kyxsxty,  & xhx^*kxysytx  \approx xhx^*kyxsytx,  \\
\end{array}\qquad \right\}\\
\left. \begin{array}[c]{rcl}
x^*hxkxysx^*ty^*  \approx x^*hxkyxsx^*ty^*,& x^*hxkxysy^*tx^*  \approx  x^*hxkyxsy^*tx^*,\\
xhx^*kxysx^*ty^*  \approx xhx^*kyxsx^*ty^*, & xhx^*kxysy^*tx^*  \approx  xhx^*kyxsy^*tx^*,
\end{array} \right\} \label{e7}
\end{align}
and their reverses constitute an identity basis for  $(\baxt_2,~^{\sharp})$.
\end{theorem}

\begin{proof}
Clearly, $(\baxt_2,~^{\sharp})$ satisfies \eqref{id: inv}, and the identities \eqref{e1}--\eqref{e7} and their reverses by Theorem~\ref{thm:baxt2}.

Note that the identities \eqref{id: inv} can be used to convert any non-empty term into some unique word. It suffices to show that each non-trivial word identity satisfied by $(\mathsf{baxt}_2,~^\sharp)$ can be deduced from \eqref{e1}--\eqref{e7} and their reverses. Note that any identity satisfied by $(\mathsf{baxt}_2,~^\sharp)$ is balanced.
Let $\Sigma$ be the set of all non-trivial balanced word identities satisfied by $(\baxt_2,~^{\sharp})$ but can not be deducible from \eqref{e1}--\eqref{e7} and their reverses. Suppose that $\Sigma\neq \emptyset$. Note that any balanced identity $\bfu\approx \bfv$ in $\Sigma$ can be written
uniquely into the form $\bfu' a\bfw \approx \bfv' b\bfw$ where $a,b$ are distinct variables and $|\bfu'|=|\bfv'|$. Choose an identity, say $\bfu\approx \bfv$, from $\Sigma$ such that when it is written as $\bfu' a\bfw\approx \bfv' b\bfw$,
the lengths of the words $\bfu'$ and $\bfv'$ are as short as possible. Since $\bfu\approx \bfv$ is balanced, we have $\con(\bfu'a)=\con(\bfv'b)$ and $\occ(x, \bfu'a)=\occ(x, \bfv'b)$ for any $x \in \con(\bfu'a)$.

Since $\con(\bfu'a)=\con(\bfv'b)$, we have $b\in \con(\bfu')$, and so $\bfu'a=\bfu_1 bc\bfu_2$ with $b\neq c$ and $b\not\in\con(\bfu_2)$.  Since $\con(\bfu'a)=\con(\bfv'b)$, we have $c\in \con(\bfv')$, and so $\bfv'=\bfv_1 c\bfv_2$ with $ c\not\in\con(\bfu_2)$. Thus $\bfu=\bfu_1 bc\bfu_2\bfw$ and $\bfv=\bfv_1 c\bfv_2b\bfw$.
On the one hand, we claim that at least one of four cases $b,b^*\in\con(\bfu_1)$, $c,c^*\in\con(\bfu_1)$, $b^*,c^*\in\con(\bfu_1)$ or $b, c\in \con(\bfu_1)$ is true. If $b, b^*, c, c^* \not\in\con(\bfu_1)$, then $\occ(b, \bfu'a)=1$, and so $\bfu[b,c]$ starts with $b$ but $\bfv[b,c]$ does not, which contradicts with (I) of Theorem~\ref{thm:baxt2}. If only $b$ occurs in  $\bfu_1$, then $\bfu[b,c]$ starts with $b^{\occ(b, \bfu'a)}c$ but $\bfv[b,c]$ does not, which contradicts with (IIa) of Theorem~\ref{thm:baxt2}. If only $b^*$ occurs in  $\bfu_1$, then $\bfu[b,c]$ starts with $(b^*)^{\occ(b^*, \bfu_1)}b$ but $\bfv[b,c]$ does not, which contradicts with (Ia) of Theorem~\ref{thm:baxt2}. If only $c$ occurs in  $\bfu_1$, then $\bfu[b,c]$ starts with $c^{\occ(c, \bfu_1)}b$ but $\bfv[b,c]$ does not, which contradicts with (IIa) of Theorem~\ref{thm:baxt2}.
If only $c^*$ occurs in $\bfu_1$, then $\bfu[b,c]$ starts with $(c^*)^{\occ(c^*, \bfu_1)}b$ but $\bfv[b,c]$ does not, which contradicts with (IIa) of Theorem~\ref{thm:baxt2}. If only $c, b^*$  occurs in $\bfu_1$, then $\bfu[b,c]$ starts with $\bfu'_1b$ with $\con(\bfu'_1)=\{c, b^*\}$ but $\bfv[b,c]$ does not starts with $\bfv'_1b$ or $\bfv'_1c^*$ with $\con(\bfv'_1)=\{c, b^*\}$ and $\occ(c, \bfv'_1)=\occ(c, \bfu'_1)$, which contradicts with (IIIa) of Theorem~\ref{thm:baxt2}. If only $c^*, b$ occurs in  $\bfu_1$, then $\bfu[b,c]$ starts with $\bfu'_1c$ with $\con(\bfu'_1)=\{c^*, b\}$ but $\bfv[b,c]$ does not starts with $\bfv'_1b^*$ or $\bfv'_1c$ with $\con(\bfv'_1)=\{c^*, b\}$ and $\occ(b, \bfv'_1)=\occ(b, \bfu'_1)$, which contradicts with (IIIa) of Theorem~\ref{thm:baxt2}.

On the other hand,
we claim that at least one of four cases $b,b^*\in\con(\bfu_2\bfw)$, $c,c^*\in\con(\bfu_2\bfw)$, $b^*,c^*\in\con(\bfu_2\bfw)$ or $b, c\in \con(\bfu_2\bfw)$ is true. If $b, b^*, c, c^* \not\in\con(\bfu_2\bfw)$, then $\bfu[b,c]$ ends with $c$ but $\bfv[b,c]$ ends with $b$, which contradicts with (I) of Theorem~\ref{thm:baxt2}. If only $b$ occurs in  $\bfu_2\bfw$, then $\bfu[b,c]$ ends with $cb^{\occ(b, \bfw)}$ but $\bfv[b,c]$ does not, which contradicts with (IIb) of Theorem~\ref{thm:baxt2}. If only $b^*$ occurs in $\bfu_2\bfw$, then $\bfu[b,c]$ ends with $c(b^*)^{\occ(b^*, \bfu_2\bfw)}$ but $\bfv[b,c]$ does not, which contradicts with (IIb) of Theorem~\ref{thm:baxt2}.
If only $c$ occurs in  $\bfu_2\bfw$, then $\bfu[b,c]$ ends with $bc^{1+\occ(c, \bfu_2\bfw)}$ but $\bfv[b,c]$ does not,  which contradicts with (IIb) of Theorem~\ref{thm:baxt2}.
If only $c^*$ occurs in  $\bfu_2\bfw$, then $\bfu[b,c]$ ends with $c(c^*)^{\occ(c^*, \bfu_2\bfw)}$ but $\bfv[b,c]$ does not, which contradicts with (Ib) of Theorem~\ref{thm:baxt2}.
If only $c, b^*$  occurs in  $\bfu_2\bfw$, then $\bfu[b,c]$ ends with $b\bfu'_2$ with $\con(\bfu'_2)=\{c, b^*\}$ but $\bfv[b,c]$ does not ends with $b\bfv'_2$ or $c^*\bfv'_2$ with $\con(\bfv'_2)=\{c, b^*\}$ and $\occ(c, \bfv'_2)=\occ(c, \bfu'_2)$, which contradicts with (IIIb) of Theorem~\ref{thm:baxt2}. If only $c^*, b$  occurs in  $\bfu_2\bfw$, then $\bfu[b,c]$ ends with $c\bfu'_2$ with $\con(\bfu'_2)=\{c^*, b\}$ but $\bfv[b,c]$ does not ends with $c\bfv'_2$ or $b^*\bfv'_2$ with $\con(\bfv'_2)=\{c^*, b\}$ and $\occ(b, \bfv'_2)=\occ(b, \bfu'_2)$, which contradicts with (IIIb) of Theorem~\ref{thm:baxt2}.
Therefore, the
identities \eqref{e1}--\eqref{e7} and their reverses can be used to convert $\bfu_1 bc\bfu_2 \bfw$ into $\bfu_1cb\bfu_2\bfw$. By repeating this process, the word $\bfu=\bfu'a\bfw=\bfu_1bc\bfu_2 \bfw$ can be converted into the word $\bfu_1c\bfu_2b \bfw$ by the
identities \eqref{e1}--\eqref{e7} and their reverses.

Clearly the identities $\bfu_1c\bfu_2b\bfw \approx \bfu'a\bfw \approx \bfv'b\bfw$ holds in $(\baxt_2, ~^{\sharp})$. Note that $|\bfu_1c\bfu_2b\bfw|=|\bfu'a\bfw|=|\bfv'b\bfw|$ and words in the identity $\bfu_1c\bfu_2b\bfw \approx \bfv'b\bfw$ have a longer common suffix than words in the identity $\bfu'a\bfw \approx \bfv'b\bfw$. Hence $\bfu_1c\bfu_2b\bfw \approx \bfv'b\bfw \notin \Sigma$ by the minimality assumption on the lengths $\bfu'$ and $\bfv'$, that is, $\bfu_1c\bfu_2b\bfw \approx \bfv'b\bfw$ can be deducible from \eqref{e1}--\eqref{e7} and their reverses. We have shown that $\bfu'a\bfw \approx \bfu_1c\bfu_2b\bfw$ can be deducible from \eqref{e1}--\eqref{e7} and their reverses. Hence $\bfu'a\bfw \approx \bfv'b\bfw$ can be deducible from \eqref{e1}--\eqref{e7} and their reverses, which contradicts with $\bfu'a\bfw \approx \bfv'b\bfw \in \Sigma$. Therefore $\Sigma=\emptyset$.
\end{proof}

Next, we show that $(\mathsf{baxt}_3,~^\sharp)$ is non-finitely based.

\begin{theorem}\label{thm: conNFB}
Suppose that $(M, \op)$ is an involution monoid satisfying the following conditions:
\begin{enumerate}[\rm(I)]
  \item for each $k \geq 2$, $(M, \op)$ satisfies the identity $\bfp_k \approx \bfq_k$ where
  \begin{align*}
\bfp_k&=x_{1}^*x_2^*\cdots x_{2k}^*\cdot xx^*\cdot x^* x_1\cdots x_{2k} x\cdot x^*x\cdot x_1^*x_3^*\cdots x_{2k-1}^*x_2^*x_4^*\cdots x_{2k}^*,\\
\bfq_k&=x_{1}^*x_2^*\cdots x_{2k}^*\cdot xx^*\cdot x x_1\cdots x_{2k} x^*\cdot x^* x\cdot  x_1^*x_3^*\cdots x_{2k-1}^*x_2^*x_4^*\cdots x_{2k}^*;
\end{align*}
  \item $\ip(\bfu)=\ip(\bfv), \fp(\bfu)=\fp(\bfv)$ and $\bfu\approx \bfv$ is balanced for any $\bfu\approx \bfv$ satisfied by $(M, \op)$;
  \item if $(M, \op)$ satisfies a non-trival identity $y^*xx^*\cdot x^*  y  x\cdot x^*x y^*\approx \bfw$, then $
  \bfw=y^*xx^*\cdot x  y  x^*\cdot x^*x y^*$;
   \item \begin{enumerate}[\rm(a)]
          \item  if $\bfu[x,y]\in \bfa y\{x,x^*,y,y^*\}^{\times}$ with $\con(\bfa)= \{x, x^*\}$, then
          $\bfv[x,y]\in \bfa' y\{x,x^*,y,y^*\}^{\times}$
           where $\bfa'$ is a permutation of $\bfa$,
          \item   if $\bfu[x,y]\in \{x, x^*, y, y^*\}^{\times} y\bfa$ with $\con(\bfa)= \{x, x^*\}$, then
          $\bfv[x,y]\in \{x,x^*,y,y^*\}^{\times}y\bfa'$
            where $\bfa'$ is a permutation of $\bfa$.
        \end{enumerate}
\end{enumerate}
Then $(M, \op)$ is non-finitely based.
\end{theorem}

The proof of Theorem~\ref{thm: conNFB} is given at the end of Lemma~\ref{lem:t=pk}.

For each $k\geq 2$, define
 \begin{align*}
\mathsf{P}_k&=x_{1}^*x_2^*\cdots x_{2k}^*\cdot xx^*\cdot x^* x_{1\pi}\cdots x_{2k\pi} x\cdot x^*x\cdot x_1^*x_3^*\cdots x_{2k-1}^*x_2^*x_4^*\cdots x_{2k}^*,\\
\mathsf{Q}_k&=x_{1}^*x_2^*\cdots x_{2k}^*\cdot xx^*\cdot x x_{1\sigma}\cdots x_{2k\sigma} x^*\cdot x^* x\cdot  x_1^*x_3^*\cdots x_{2k-1}^*x_2^*x_4^*\cdots x_{2k}^*,
\end{align*}
where $\pi, \sigma$ are any permutations on $\{1, 2, \dots, 2k\}$.

\begin{lemma}\label{pkqk}
Let $(M, \op)$ be an involution monoid satisfying conditions (II) and (III)
in Theorem \ref{thm: conNFB}. Suppose that $\bfp_k \approx \bfw$ is any word identity satisfied by $(M, \op)$ such that $\bfp_k \in \mathsf{P}_k$ . Then $\bfw \in \mathsf{P}_k \cup \mathsf{Q}_k$.
\end{lemma}

\begin{proof}
Note that $\bfp_k[x_i, x]=x_i^*xx^*\cdot x^* x_i x\cdot x^*xx_i^*$ for $i=1,2,\cdots, 2k$. It follows from the condition (II) of Theorem \ref{thm: conNFB} that $\bfw[x_i, x]=x_i^*xx^*\cdot x^* x_i x\cdot x^*xx_i^*$ or $\bfw[x_i, x]=x_i^*xx^*\cdot x x_i x^*\cdot x^*xx_i^*$. Then it follows from (III) Theorem \ref{thm: conNFB} that $\bfw \in \mathsf{P}_k \cup \mathsf{Q}_k$.
\end{proof}

A word identity $\bfu \approx \bfv$ is \textit{k-limited} if $\con(\overline{\bfu\bfv}) \leq k$. For
any involution semigroup $(S, ~^*)$, let $\mathsf{id_k}(S, ~^*)$ denote the set of all $k$-limited word
identities of $(S, ~^*)$.

\begin{lemma}\label{lem:t=pk}
Suppose that $(M, \op)$ is an involution monoid satisfying conditions (II)--(IV)
in Theorem \ref{thm: conNFB}. Let $\bfs \approx \bft$ be an identity which is directly deducible from some identity in $\mathsf{id}_{2k}(M, \op)$ with $\lfloor\bfs\rfloor \in \mathsf{P}_k$. Then $\lfloor\bft\rfloor \in \mathsf{Q}_k$.
\end{lemma}

\begin{proof}
Let $\bfu \approx \bfv$ be a word identity in $\mathsf{id}_{2k}(S, \op)$ from which the identity $\bfs \approx \bft$ is
directly deducible. There is a substitution $\varphi: \mathcal{X} \rightarrow \mathsf{T}(\mathcal{X})$ such that $\varphi(\bfu)$ is a subterm
of $\bfs$, and replacing this particular subterm $\varphi(\bfu)$ of $\bfs$ with $\varphi(\bfv)$ results in $\bft$. By Remark \ref{rem:factor},
either $\lfloor\varphi(\bfu)\rfloor$ or $\lfloor(\varphi(\bfu))^*\rfloor$ is a factor of
$\lfloor \bfs\rfloor$. It suffices to consider the former case since the
latter is similar. Hence there exist words $\bfa, \bfb \in (\mathcal{X}\cup \mathcal{X}^*)^{+}$ such that
 $\lfloor\bfs\rfloor = \bfa \lfloor \varphi(\bfu)\rfloor \bfb$.
Since $\bft$ is obtained by replacing $\varphi(\bfu)$ in $\bfs$ with $\varphi(\bfv)$, it follows that
$\lfloor\bft \rfloor= \bfa \lfloor \varphi(\bfv)\rfloor \bfb$. Since $\lfloor\bfs\rfloor \in \mathsf{P}_k$, it follows from Lemma \ref{pkqk} that $\lfloor\bft \rfloor \in \mathsf{P}_k\cup \mathsf{Q}_k$. Working
toward a contradiction, suppose that $\lfloor\bft \rfloor \in \mathsf{Q}_k$. Since
$\lfloor \bfs\rfloor $ and $\lfloor \bft\rfloor $ share the same prefix $\bfa$ and the same suffix $\bfb$ and $\ip(\lfloor\varphi(\bfu)\rfloor)=\ip(\lfloor\varphi(\bfv)\rfloor)$ and $\fp(\lfloor\varphi(\bfu)\rfloor)=\fp(\lfloor\varphi(\bfv)\rfloor)$, $\lfloor\varphi(\bfu)\rfloor$ contains the factor
\[
x^*\cdot x^* x_{1\pi}\cdots x_{2k\pi} x\cdot x^*
\]
and $\lfloor \bfv\varphi\rfloor$ contains the factor
\[
x^*\cdot x x_{1\pi}\cdots x_{2k\pi} x^*\cdot x^*.
\]
It follows from (IV) of Theorem \ref{thm: conNFB} that $\overleftarrow{\occ}_{x_{i\pi}}(x_{i\pi}^*, \bfu\varphi)\neq 0$ and $\overrightarrow{\occ}_{x_{i\pi}}(x_{i\pi}^*, \varphi(\bfu))\neq 0$ for each $i\in\{1, 2, \dots, 2k\}$, and so $\bfa, \bfb=\emptyset$.
Since $\bfu\approx \bfv$ is $2k$-limited and $\lfloor\varphi(\bfu)\rfloor\approx \lfloor\varphi(\bfv)\rfloor$ is $2k+1$-limited, there
exists a variable $s \in \con(\bfu)$ such that $\lfloor\varphi(s)\rfloor$ contains one of the following factors:
\[
x^*x_{1\pi} ,\,\, x_{1\pi}x_{2\pi}, \,\, x_{2\pi}x_{3\pi},\,\, \cdots, \,\, x_{2k-1\pi}x_{2k\pi},\,\, x_{2k\pi}x.
\]
Suppose that there exists a variable $s\in\con(\bfu)$ such that $\lfloor \varphi(s) \rfloor$ contains $x_{2k\pi}x$ as a factor. Then $\lfloor\bft\rfloor$ also contains $x_{2k\pi}x$ as a factor, which contradicts with $\lfloor\bft \rfloor \in\mathsf{Q}_k$. Similarly, there does not exist $s\in\con(\bfu)$ such that $\lfloor \varphi(s) \rfloor$ contains $x^*x_{1\pi}$ as a factor. Suppose that there exists a variable $s \in \con(\bfu)$ such that $\lfloor\varphi(s)\rfloor$ contains $x_{i\pi}x_{(i+1)\pi}$ as a factor for some $i=1, 2, \dots, 2k-1$.
Clearly, $\occ(s, \bfu)=1$. Then $x_{i\pi}x_{(i+1)\pi}$ is a factor of
$\lfloor \varphi(\bfv) \rfloor$ since $\bfu\approx \bfv$ is balanced.
It follows from the definition of $\mathsf{P}_k$ that $\bfu=\bfu_1s\bfu_2$, and $s^*$ can not occur in $\con(\bfu_1)$ and $\con(\bfu_2)$ simultaneously. Suppose that $s^*\not\in\con(\bfu_2)$. Note that $xx^*x$ is a factor of $\lfloor\bfs\rfloor$ occurring after $\lfloor \varphi(s)\rfloor$. Then there exist at most three variables in $\con(\bfu_2)$, say $s_1, s_2, s_3$, such that $\lfloor\varphi(s_1s_2s_3)\rfloor$ contains $xx^*x$ as a factor in $\lfloor\bfs\rfloor$. It follows from conditions (II) and (IV) of Theorem \ref{thm: conNFB} that $\bfv=\bfv_1s\bfv_2$ satisfying $s, s^* \not \in \con(\bfv_2), \occ(s_1,\bfv_2)=\occ(s_1,\bfu_2), \occ(s_2, \bfv_2)=\occ(s_1,\bfu_2)$ and $\occ(s_3, \bfv_2)=\occ(s_3,\bfu_2)$, and so $xx^*x$ occurs after $x_{i\pi}x_{i+1\pi}$ in $\lfloor\bft\rfloor$,
but this contradicts with $\bfa\lfloor \varphi(\bfv)\rfloor \bfb \in \mathsf{Q}_k$.
Consequently, $\bfa\lfloor \varphi(\bfv)\rfloor \bfb \in \mathsf{P}_k$.
\end{proof}

\begin{proof}[\bf{Proof of Theorem \ref{thm: conNFB}}]
Let $(M,~^*)$ be any involution monoid that satisfies conditions (I)--(IV) in Theorem \ref{thm: conNFB}. Then there exists some set $\Sigma$ of word identities such that \eqref{id: inv}$\cup \Sigma$
is an identity basis for $(M,~^*)$. Working toward a contradiction, suppose that $(M,~^*)$
is finitely based. Then there exists a finite subset  $\Sigma_{\mathsf{fin}}$ of $\Sigma$  such that all
identities of $(M,~^*)$ are deducible from \eqref{id: inv}$\cup\Sigma_{\mathsf{fin}}$. Hence there exists some fixed integer
$k$ such that  $\Sigma_{\mathsf{fin}} \subseteq\Sigma\cap \mathsf{id}_{2k} (M,~^*)$. By (I), the involution monoid $(M,~^*)$ satisfies some
word identity $\bfp_k \approx \bfq_k$ with $\bfp_k \in \mathsf{P}_k$ and $\bfq_k \in \mathsf{Q}_k$. Therefore there
exists some sequence
\[
\bfp_k = \bfs_1,\bfs_2,\cdots,\bfs_m = \bfq_k
\]
of terms such that each identity $\bfs_i\approx \bfs_{i+1}$ is directly deducible from some identity
$\bfu_i \approx \bfv_i \in $ \eqref{id: inv}$\cup \Sigma_{\mathsf{fin}}$. The equality
$\bfs_1 = \bfp_k \in \mathsf{P}_k$ holds. If $\bfs_i  \in \mathsf{P}_k$ for some $i\geq 1$, then there are two cases depending on whether the
identity $\bfu_i \approx \bfv_i $ is from \eqref{id: inv} or $ \Sigma_{\mathsf{fin}}$. If $\bfu_i \approx \bfv_i $ is from \eqref{id: inv}, then
$\bfs_i = \bfs_{i+1}$ by Remark \ref{rem:s=t}, whence $\bfs_{i+1}\in \mathsf{P}_k$. If $\bfu_i \approx \bfv_i $ is from $\Sigma_{\mathsf{fin}}$, then $\bfs_{i+1}\in \mathsf{P}_k$  by
Lemma \ref{lem:t=pk}. Therefore $\bfs_{i+1}\in \mathsf{P}_k$ in any case, whence by induction,
$\bfs_{i}\in \mathsf{P}_k$ for all $i$. But this implies the contradiction $\bfq_k = \bfs_m \in \mathsf{P}_k$. Consequently, the involution monoid $(M,~^*)$ is non-finitely based.
\end{proof}

\begin{theorem}
The involution monoid $(\mathsf{baxt}_3,~^\sharp)$ is non-finitely based.
\end{theorem}

\begin{proof}
It follows from Theorem \ref{thm:baxt3} that $(\mathsf{baxt}_3,~^\sharp)$ satisfies all of  conditions of Theorem \ref{thm: conNFB}. Therefore $(\mathsf{baxt}_3,~^\sharp)$ is non-finitely based.
\end{proof}

Now, we consider the finite basis problem for $(\baxt_n,~^{\sharp})$ with $n\geq 4$.

\begin{theorem}\label{thm:baxt4's basis}
The identities \eqref{id: inv} and
\begin{align}
xhyk\,xy\,sxty  \approx xhyk\,yx\,sxty, \quad & xhyk\,xy\,sytx  \approx xhyk\,yx\,sytx,\tag{5.3}\label{e}
\end{align}
constitute an identity basis for  $(\baxt_n,~^{\sharp})$ for $n\geq 4$.
\end{theorem}

\begin{proof}
By Theorem \ref{thm:same-var}, we only need to show that the result holds for $(\baxt_4,~^{\sharp})$.
Clearly, $(\baxt_4,~^{\sharp})$ satisfies the identities \eqref{id: inv}, and \eqref{e} by Theorem~\ref{thm:baxt4+}.

Note that the identities \eqref{id: inv} can be used to convert any non-empty term into some unique word. It suffices to show that each non-trivial word identity satisfied  by $(\mathsf{baxt}_4,~^\sharp)$ can be derived from \eqref{e}.
Note that any identity satisfied by $(\mathsf{baxt}_4,~^\sharp)$ is balanced.
Let $\Sigma$ be the set of all non-trivial balanced word identities satisfied by $(\baxt_4,~^{\sharp})$ but can not be deducible from \eqref{e}. Suppose that $\Sigma\neq \emptyset$. Note that any balanced identity $\bfu\approx \bfv$ in $\Sigma$ can be written
uniquely into the form $\bfu' a\bfw \approx \bfv' b\bfw$ where $a,b$ are distinct variables and $|\bfu'|=|\bfv'|$. Choose an identity, say $\bfu\approx \bfv$, from $\Sigma$ such that when it is written as $\bfu' a\bfw\approx \bfv' b\bfw$,
the lengths of the words $\bfu'$ and $\bfv'$ are as short as possible. Since $\bfu\approx \bfv$ is balanced, we have $\con(\bfu'a)=\con(\bfv'b)$ and $\occ(x, \bfu'a)=\occ(x, \bfv'b)$ for any $x \in \con(\bfu'a)$.

Since $\con(\bfu'a)=\con(\bfv'b)$, we have $b\in \con(\bfu')$, and so $\bfu'a=\bfu_1 bc\bfu_2$ with $b\neq c, b\not\in\con(\bfu_2)$.  Since $\con(\bfu'a)=\con(\bfv'b)$, we have $c\in \con(\bfv')$, and so $\bfv'=\bfv_1 c\bfv_2$ with $ c\not\in\con(\bfu_2)$. Thus $\bfu=\bfu_1 bc\bfu_2\bfw$ and $\bfv=\bfv_1 c\bfv_2b\bfw$.
We claim that $b, c\in \con(\bfu_1)$. If $b, c \not\in\con(\bfu_1)$, then $\occ(b, \bfu'a)=1$, and so $\overleftarrow{\occ}_c(b, \bfu)=1$ but $\overleftarrow{\occ}_c(b, \bfv)=0$, which contradicts with Theorem~\ref{thm:baxt4+}. If only $b$ occurs in  $\bfu_1$, then $\overleftarrow{\occ}_c(b, \bfu)>\overleftarrow{\occ}_c(b, \bfv)$, which contradicts with Theorem~\ref{thm:baxt4+}. If only $c$ occurs in  $\bfu_1$, then $\overleftarrow{\occ}_b(c, \bfu)<\overleftarrow{\occ}_b(c, \bfv)$, which contradicts with Theorem~\ref{thm:baxt4+}.
We also claim that $b, c\in \con(\bfu_2\bfw)$. If $b, c\not\in\con(\bfu_2\bfw)$, then $\overrightarrow{\occ}_b(c, \bfu)=1$ but $\overrightarrow{\occ}_b(c, \bfv)=0$, which contradicts with Theorem~\ref{thm:baxt4+}. If only $b$ occurs in  $\bfu_2\bfw$, then $\overrightarrow{\occ}_c(b, \bfu)<\overrightarrow{\occ}_c(b, \bfv)$, which contradicts with Theorem~\ref{thm:baxt4+}.
If only $c$ occurs in  $\bfu_2\bfw$, then $\overrightarrow{\occ}_b(c, \bfu)>\overrightarrow{\occ}_b(c, \bfv)$,  which contradicts with Theorem~\ref{thm:baxt4+}.
As such, we can deduce the word $\bfu_1cb\bfu_2\bfw$, by applying the
identities \eqref{e}. By repeating this process, the word $\bfu_1bc\bfu_2 \bfw$ can be converted into the word $\bfu_1c\bfu_2b \bfw$.

Clearly the identities $\bfu_1c\bfu_2b\bfw \approx \bfu'a\bfw \approx \bfv'b\bfw$ holds in $(\baxt_4, ~^{\sharp})$. Note that $|\bfu_1c\bfu_2b\bfw|=|\bfu'a\bfw|=|\bfv'b\bfw|$ and words in the identity $\bfu_1c\bfu_2b\bfw \approx \bfv'b\bfw$ have a longer common suffix than words in the identity $\bfu'a\bfw \approx \bfv'b\bfw$. Hence $\bfu_1c\bfu_2b\bfw \approx \bfv'b\bfw \notin \Sigma$ by the minimality assumption on the lengths of $\bfu'$ and $\bfv'$, that is, $\bfu_1c\bfu_2b\bfw \approx \bfv'b\bfw$ can be deducible from \eqref{e}. We have shown that $\bfu'a\bfw \approx \bfu_1c\bfu_2b\bfw$ can be deducible from \eqref{e}. Hence $\bfu'a\bfw \approx \bfv'b\bfw$ can be deducible from \eqref{e}, which contradicts with $\bfu'a\bfw \approx \bfv'b\bfw \in \Sigma$. Therefore $\Sigma=\emptyset$.
\end{proof}

In the end, we consider the number of subvarieties of $\var(\baxt_n,~^{\sharp})$ for each $n\geq 2$.
Recall that a word $\bfu$ is an \textit{isoterm} for an involution monoid if it does not satisfy any non-trivial word
identity of the form $\bfu \approx \bfv$.

\begin{lemma}[{\cite[Theorem~3.6]{GZL-variety}}]\label{lem:subvariety}
Let $(M,~^*)$ be any involution monoid with isoterms $xx^*yy^*$ and $xyy^*x^*$.
Then the variety $\var(M,~^*)$ contains continuum many subvarieties.
\end{lemma}

\begin{theorem}\label{th:subvarieties}
The variety $\var(\baxt_n,~^{\sharp})$ for each finite $n\geq 2$ contains continuum many subvarieties.
\end{theorem}

\begin{proof}
By Theorems~\ref{thm:baxt2}--\ref{thm:baxt4+}, it is routine to show that both the words $xx^*yy^*$ and $xyy^*x^*$ are isoterms for  $(\baxt_n,~^{\sharp})$ for each finite $n\geq 2$. Now the results follows from Lemma~\ref{lem:subvariety}.
\end{proof}

\section{Recognizing identities of $(\mathsf{baxt}_n,~^\sharp)$ in polynomial time}\label{sec:CHECK-ID}%

In this section, it is shown that the identity checking problem of $(\mathsf{baxt}_n,~^\sharp)$ for each finite $n$
belongs to the complexity class $\mathsf{P}$.

Let $\bfu\in (\mathcal{X}\cup \mathcal{X}^*)^{+}$.
If $x, x^* \in \con (\bfu)$, then $\{x, x^*\}$ is called a \textit{mixed pair} of $\mathbf{u}$.
Denote by $\mathsf{pre}(\bfu)$ the longest prefix of $\bfu$ containing only one variable, that is, $|\con(\mathsf{pre}(\bfu))|=1$, and $\mathsf{suf}(\bfu)$ the longest suffix of $\bfu$ containing only one variable, that is, $|\con(\mathsf{suf}(\bfu))|=1$;
and denote by $\mathsf{pren}(\bfu)$ the longest prefix of $\bfu$ which does not contain any mixed pair, and $\mathsf{sufn}(\bfu)$ the longest suffix of $\bfu$ which does not contain any mixed pair.

\begin{theorem}\label{thm:checkbaxt}
The decision problem ${\textsc{Check-Id}}(\mathsf{baxt}_n,~^\sharp)$ for each finite $n$
belongs to the complexity class $\mathsf{P}$.
\end{theorem}

\begin{proof}
For decision problem ${\textsc{Check-Id}}(\mathsf{baxt}_1,~^\sharp)$, given any word identity $\bfu\approx \bfv$, it suffices to show that  whether $\con(\overline{\bfu})= \con(\overline{\bfv})$ and $\occ(x, \bfu)+\occ(x^*, \bfu)=\occ(x, \bfv)+\occ(x^*, \bfv)$ for any $x\in \con(\overline{\bfu})$. Clearly, these can be completed in polynomial time.

For decision problem ${\textsc{Check-Id}}(\mathsf{baxt}_2,~^\sharp)$, given any word identity $\bfu\approx \bfv$, it suffices to show that one can check whether the identity $\bfu \approx \bfv$ is balanced and the words $\bfu$ and  $\bfv$ satisfy conditions of Theorem \ref{thm:baxt2} in polynomial time.

To check whether $\bfu \approx \bfv$ is balanced, it suffices to check whether $\con(\bfu) =\con(\bfv)$ and $\occ(x, \bfu)=\occ(x, \bfv)$ for any $x\in \con(\bfu)$. Clearly, these can be completed in polynomial time.

To check  conditions (I)--(II) of Theorem \ref{thm:baxt2}, it suffices to check that, for any $x, y \in \con(\bfu)$ with $x, x^*\neq y$, whether $\mathsf{pre}(\bfu[x,y])=\mathsf{pre}(\bfv[x,y])$, and the variable following immediately $\mathsf{pre}(\bfu[x,y])$ is the same as the variable following immediately $\mathsf{pre}(\bfv[x,y])$; and to check whether $\mathsf{suf}(\bfu[x,y])=\mathsf{suf}(\bfv[x,y])$, and the variable preceding immediately $\mathsf{suf}(\bfu[x,y])$ is the same as the variable preceding immediately $\mathsf{suf}(\bfv[x,y])$.  There are at most $\binom{|\con(\bfu)|}{2}$ pairs of such $x, y$. Clearly, these can be completed in polynomial time.

To check condition (III) of Theorem \ref{thm:baxt2}, it suffices to check that, for any $x, y \in \con(\bfu)$ with $x, x^*\neq y$,
whether $\con(\mathsf{pren}(\bfu[x,y]))=\con(\mathsf{pren}(\bfv[x,y]))$ and $\occ(z, \mathsf{pren}(\bfu[x,y]))=\occ(z, \mathsf{pren}(\bfv[x,y]))$  for any $z\in\con(\mathsf{pren}(\bfu[x,y]))$; and whether $\con(\mathsf{sufn}(\bfu[x,y]))=\con(\mathsf{sufn}(\bfv[x,y]))$ and $\occ(z, \mathsf{sufn}(\bfu[x,y]))=\occ(z, \mathsf{sufn}(\bfv[x,y]))$  for any $z\in\con(\mathsf{sufn}(\bfu[x,y]))$.  There are at most $\binom{|\con(\bfu)|}{2}$ pairs of such $x, y$.
Clearly, these can be completed in polynomial time.

Therefore the decision problem ${\textsc{Check-Id}}(\mathsf{baxt}_2,~^\sharp)$
belongs to the complexity class $\mathsf{P}$.

For decision problem ${\textsc{Check-Id}}(\mathsf{baxt}_3,~^\sharp)$, given any word identity $\bfu\approx \bfv$, by Theorem \ref{thm:baxt2} and the above arguments, we only need to show that one can check whether the words $\bfu$ and  $\bfv$ satisfy conditions (III)--(V) of Theorem \ref{thm:baxt3} in polynomial time.

To check condition (III) of Theorem \ref{thm:baxt3}, it suffices to check that, for any $x, y \in \con(\bfu)$ with $x, x^*\neq y$,
whether $\con(\mathsf{pren}(\bfu[x,y]))=\con(\mathsf{pren}(\bfv[x,y]))$, $\occ(z, \mathsf{pren}(\bfu[x,y]))=\occ(z, \mathsf{pren}(\bfv[x,y]))$  for any $z\in\con(\mathsf{pren}(\bfu[x,y]))$ and the variable following immediately $\mathsf{pren}(\bfu[x,y])$ is the same as the variable following immediately $\mathsf{pren}(\bfv[x,y])$; and to check whether $\con(\mathsf{sufn}(\bfu[x,y]))=\con(\mathsf{sufn}(\bfv[x,y]))$, $\occ(z, \mathsf{sufn}(\bfu[x,y]))=\occ(z, \mathsf{sufn}(\bfv[x,y]))$  for any $z\in\con(\mathsf{sufn}(\bfu[x,y]))$ and the variable preceding immediately $\mathsf{pren}(\bfu[x,y])$ is the same as the variable preceding immediately $\mathsf{pren}(\bfv[x,y])$.  There are at most $\binom{|\con(\bfu)|}{2}$ pairs of such $x, y$.
Clearly, these can be completed in polynomial time.

To check condition (IV) of Theorem \ref{thm:baxt3}, it suffices to check that, for any $x, y\in \con(\bfu)$ with $x, x^*\neq y$,
whether $\overleftarrow{\occ}_y(x, \bfu)+\overleftarrow{\occ}_y(x^*, \bfu)= \overleftarrow{\occ}_y(x, \bfv)+\overleftarrow{\occ}_y(x^*, \bfv)$ and $\overrightarrow{\occ}_y(x, \bfu)+\overrightarrow{\occ}_y(x^*, \bfu)= \overrightarrow{\occ}_y(x, \bfv)+\overrightarrow{\occ}_y(x^*, \bfv)$. There are at most $\binom{|\con(\bfu)|}{2}$ pairs of such $x, y$.  Clearly, these can be completed in polynomial time.

To check condition (V) of Theorem \ref{thm:baxt3}, it suffices to check that, for any $x, y\in \con(\bfu)$ with $x, x^*\neq y$ satisfying $\overleftarrow{\occ}_y(y^*, \bfu)=0$  whether $\overleftarrow{\occ}_y(x, \bfu)=\overleftarrow{\occ}_y(x, \bfv)$ and $\overleftarrow{\occ}_y(x^*, \bfv)=\overleftarrow{\occ}_y(x^*, \bfv)$; and for any $x, x^*\neq y\in \con(\bfu)$ satisfying $\overrightarrow{\occ}_y(y^*, \bfu)=0$, whether $\overrightarrow{\occ}_y(x, \bfu)=\overrightarrow{\occ}_y(x, \bfv)$ and $ \overrightarrow{\occ}_y(x^*, \bfv)=\overrightarrow{\occ}_y(x^*, \bfv)$. There are at most $\binom{|\con(\bfu)|}{2}$ pairs of such $x, y$. Clearly, these can be completed in polynomial time.

Therefore the decision problem ${\textsc{Check-Id}}(\mathsf{baxt}_3,~^\sharp)$
belongs to the complexity class $\mathsf{P}$.

For decision problem ${\textsc{Check-Id}}(\mathsf{baxt}_n,~^\sharp)$ with $n\geq 4$,  it suffices to check that, for any $x,y\in \con(\bfu)$,
whether $\overleftarrow{\occ}_x(y, \bfu)=\overleftarrow{\occ}_x(y, \bfv)$ and $\overrightarrow{\occ}_x(y, \bfu)=\overrightarrow{\occ}_x(y, \bfv)$. There are at most $\binom{|\con(\bfu)|}{2}$ pairs of such $x, y$. Clearly, these can be completed in polynomial time. Therefore the decision problem ${\textsc{Check-Id}}(\mathsf{baxt}_n,~^\sharp)$ with $n\geq 4$
belongs to the complexity class $\mathsf{P}$.
\end{proof}

\end{sloppypar}
\end{document}